\documentclass{amsart}
\usepackage{amsmath}
\usepackage{amsfonts}
\usepackage{amssymb}
\catcode `\@=11
\def\numberbysection{\@addtoreset{equation}{section}
         \renewcommand{\theequation}{\thesection.\arabic{equation}}}
\numberbysection
\def\subsubsection{\@startsection{subsubsection}{3}%
  \normalparindent{.5\linespacing\@plus.7\linespacing}{-.5em}%
  {\normalfont\bfseries}}
\def\a{\alpha}

\def\g{\gamma}

\def\eps{\epsilon}
\def\l{\lambda}
\def\s{\sigma}

\def\gen{\Lambda}

\def\Zz{\mathbb{Z}/2\mathbb {Z}}
\def\Zn{\mathbb{Z}/n\mathbb{Z}}
\def\Znn{\mathbb{Z}/(n+1)\mathbb{Z}}
\def\znine{\mathbb{Z}/9\mathbb{Z}}
\def\n{\{1,\dots,n\}}
\def\nn{\nonumber}

\def\Z{\mathbb{Z}}
\def\Zzn{\mathbb{Z}/(2n\mathbb{Z})}

\def\codim{\rm{codim}}

\def\dg{D(k[G])}
\def\dh{D[k[G])}

\def\cQ{\mathcal{Q}}
\def\cbQ{\bar{\mathcal{Q}}}

\def\Q{Q}
\def\M{M}

\def\credo{orbifold mirror philosophy }
\def\Credo{Orbifold mirror philosophy }
\def\cred{orbifold mirror philosophy}

\newtheorem{thm}{Theorem}[section]
\newtheorem{lm}{Lemma}[section]
\newtheorem{prop}{Proposition}[section]
\newtheorem{crl}{Corollary}[section]

\newtheorem{dfprop}[prop]{Definition-Proposition}
\newtheorem{df}{Definition}[section]
\newtheorem{rk}{Remark}[section]

\newtheorem{conj}{Conjecture}[section]
\newtheorem{phil}{Philosophy}[section]

\begin{document}

\title[Singularities, $G$--Frobenius algebras and Mirror Symmetry]
{Singularities with Symmetries,
Orbifold Frobenius Algebras and Mirror Symmetry}

\author[Ralph M. Kaufmann]{Ralph M. Kaufmann}
\email{ralphk@mpim-bonn.mpg.de}

\address{Oklahoma State University, Department of Mathematics, MS 401,
Stillwater OK 74078-1058, USA }

\begin{abstract}
Previously, we introduced a duality transformation
for Euler $G$--Frobenius
algebras. Using this transformation,
we prove that the simple $A,D,E$ singularities and
Pham  singularities of coprime powers are mirror self--dual
where the mirror duality is implemented by orbifolding with respect to
the symmetry group generated by the grading operator and dualizing.
We furthermore calculate
orbifolds and duals to other $G$--Frobenius algebras which relate
different $G$--Frobenius algebras for singularities. In particular, using
orbifolding and the duality transformation we
provide a mirror pairs for the  simple boundary singularities $B_n$ and $F_4$.
Lastly,  we relate our constructions to $r$ spin--curves, classical
singularity theory and foldings of Dynkin diagrams.

\end{abstract}
\maketitle


\section*{Introduction}
In \cite{orb} we introduced a duality transformation for Euler
$G$--Frobenius algebras which are graded Frobenius algebras whose
grading operator is realized by the action of a central element.
Using this transformation in the setting of isolated singularities
with symmetries, we prove that the simple  singularities $A,D,E$
and certain Pham singularities are mirror self--dual where the
mirror duality is implemented by orbifolding with respect to the
symmetry group generated by the grading operator and dualizing. In
particular the invariants of the orbifold are $A_1$ while the
invariants of the dual are  the simple singularity of type $A,D,E$
one started out with. Thus orbifolding and dualizing provides a
mirror dual pair to the pair $(W,A_1)$ which is naturally
associated to $W$, for $W$ one of the simple singularities
$A_n,D_n,E_6,E_7,E_8$. We also show that the same holds true for
Pham singularities of co--prime powers.

 Furthermore we calculate
orbifolds and duals to other $G$--Frobenius algebras which relate
different $G$--Frobenius algebras for singularities to each other.
We thereby provide more mirror pairs notably mirror pairs for the
simple boundary singularities. In particular
$((B_n,I_2(4)),(I_2(4),B_n))$ is obtained by orbifolding and
dualizing either $A_{2n-1}$ or $D_{n+1}$ by $\Zz$ and $\Zn$. And
$((F_4,I_2(4)),(I_2(4),F_4))$ obtained by orbifolding $E_6$ and
dualizing with respect to $\Zz$ and $\Z/3\Z \times \Zz$.

The invariants of the $G$--Frobenius algebras based on the
singularities with symmetries are related to the singularities
considered on the orbifold of $\mathbb{C}^n$ with respect to the
symmetry group while
the duals also conjecturally
play a role in the analogs of $r$ spin--curves built on
quasi--homogenous polynomials of which special types have been
studied by \cite{FJR}. For the exact formulation of these conjectures
we refer to \S \ref{spinparagraph}.

 The operation of dualizing as defined in \cite{orb}
 was inspired by the representation theory of $N=2$
super--conformal field theory applied to orbifold Landau--Ginzburg models
\cite{IV}.
Although the background is very elaborate and involves many highly complicated
concepts, in the special case we are considering all can be stated in terms
of $G$--Frobenius algebras or $\dg$ modules and algebras, where $\dg$ stands for the
Drinfel'd double of the group ring $k[G]$. $\dg$ modules
are a special type of
$G$--graded $G$--modules, namely those, where the $G$--action acts by
conjugation on the $G$--grading, cf.\ \cite{Mo,disc,JKK}.

We will first review the background for this operation and then
comment on its realization on the level of Euler $G$--Frobenius algebras.
The reader not inclined to read about physics can thus skip the following
two paragraphs and further comments about
physics which can be considered as motivation
and continue to the purely algebraic part of the paper.

A so--called $(2,2)$ super--conformal field theory has an $N=2$ super--conformal
symmetry for both the left and the right movers. This implies that there are
four finite rings which are closed under the naive operator product.
These rings are called $(c,c),(a,c),(a,a)$ and $(c,a)$ respectively.
In terms of the representation theory these rings are given by fields
which are annihilated by certain operators or
equivalently satisfy certain constraints for their eigenvalues with respect to
the operators $J_0, \bar J_0,L_0, \bar L_0$ of the two $N=2$ super--conformal
algebras, which are usually called
$q,\bar q, h$ and $\bar h$ respectively. The left $c$ or $a$ stands for left
chiral or anti--chiral and the
letter $a$ or $c$ on the right for right chiral or right anti--chiral.
An element $|\phi\rangle$ is left chiral
if  $G^+_{-1/2}| \phi \rangle=0$ or equivalently $h=\frac{q}{2}$.
It is called left anti--chiral if
$G^-_{-1/2}| \phi \rangle=0$ or equivalently $h=-\frac{q}{2}$.
Right chiral means that $\bar G^+_{-1/2}| \phi \rangle=0$ or equivalently
$\bar h=\frac{\bar q}{2}$ and finally
right anti--chiral means
that $\bar G^-_{1/2}| \phi \rangle=0$ or equivalently $\bar h=-\frac{\bar q}{2}$.
It turns out the rings $(a,a)$ and $(c,a)$ can be recovered from
$(c,c)$ and $(a,c)$ by charge conjugation. Thus one confines oneself to
study the latter two rings. Mirror symmetry as it was originally conceived
in physics was an operation which takes one conformal field theory $T$ and
produces another conformal field theory $\check T$ such that the
$(c,c)$ ring of $T$ is isomorphic to the $(a,c)$ ring of $\check T$ and
vice versa.

One special type of $N=2$ theory is given by the so--called Landau--Ginzburg
theory which is the conformally invariant fixed point of the  Lagrangian
$$
\mathcal{L}= \int   K(X,\bar X) d^2 z d^4\theta+\int f(z_i)
+ \text{ complex conjugate } d^2z d^2 \theta.
$$
where $f$ is a quasi--homogenous function of fractional degree $q_i$ for
$z_i$.
This model leads to a trivial $(a,c)$ ring and a $(c,c)$ ring which is
given by $\mathbb C[\mathbf{z}]/J_f$ where $J_f = (f_{z_i})$ is the Jacobian
ideal. Moreover the bi--degree $(q,\bar q)$ for $z_i$ is given by $(q_i,q_i)$.

The above considerations are the starting point for a purely algebraic
consideration. If the function $f$ above has an isolated
singularity at zero, the situation is one that has
 been studied for a long time by
mathematicians. The $(c,c)$ ring is in this case just the local
or Milnor ring of the singularity. The only unusual thing is the
bi--grading instead of the grading, but in fact the bi--grading is
just a diagonal grading obtained from the usual grading in singularity
theory and it contains no additional information. It will however play an
important role later on.

In the setup above, the quasi--homogeneity of the function $f$ allows
one to consider it as a function on a weighted projective space.
In the case that the polynomial describes a Calabi--Yau hypersurface
the claim these two geometries (singularity/Calabi-Yau) should give the
same Frobenius manifolds of field theories
is the famous Landau--Ginzburg/Calabi-Yau correspondence.
In doing so one is naturally considering the quotient of the theory
by a finite symmetry group. In general one can consider a group
$G\subset GL(\mathbb{C},n)$ which leaves $f(z_1,\dots, z_n)$ invariant
and consider the resulting orbifold. This particular situation and
the general setup of global orbifolds was analyzed in \cite{orb}.
It turns out that the algebraic object one is dealing with
is an extension of the Milnor ring, which by itself is a Frobenius algebra,
to a $G$--Frobenius algebra in the sense of \cite{orb}.
A $G$--Frobenius algebra has a $G$--action and
the invariants of this $G$--Frobenius are expected to form a Frobenius algebra.
These will be bi--graded in a natural way.
In physics terms this algebra of invariants is the $(c,c)$ ring
of the orbifold model. Now again appealing to physics, the orbifold
theory should also have an $(a,c)$ ring. This ring is what is computed
by the duality transformation we gave in \cite{orb}. To be precise, the ring
$(a,c)$ will be equal to the $G$--invariants of the dual $\dg$ model.
In order to define the full dual it is necessary for the group of symmetries
to contain the symmetry provided by the exponential grading operator
$J=diag(\exp (2\pi i q_1),\dots \exp(2\pi i q_n)$.

We called the transformation a mirror transformation, since as we show below,
the orbifold of the simple singularities of type $A,D,E$ by the symmetry
group generated by $J$ has a trivial
$(c,c)$ and an $(a,c)$ ring that is isomorphic to the Milnor ring of the
singularity and hence is mirror dual to the original Milnor ring.
Thus for these singularities the operation of orbifolding
and taking the invariants of the dual (i.e. the $(a,c)$ ring)
implements mirror symmetry. If one would like to phrase mirror symmetry
in terms of $A$--models and $B$--models, the Landau--Ginzburg model
is a $B$--model. In mathematical terms the $B$--Model is the Milnor ring with
the diagonal bi--grading $(q,q)$. The corresponding mirror model is an
$A$--model (not to be confused with the $A$--type singularity)
which would be given by the Milnor ring but with a grading
of $(-q,q)$. This would be a ``Landau--Ginzburg $A$--model''.


 In \cite{orb}, we have made the case that for global orbifolds it
is not enough to consider just the invariants of the $G$--action
of the $G$--Frobenius algebra, but instead one needs to
consider the whole $G$ Frobenius.
The fruitfulness of this point of view can be seen for instance
in its application to symmetric products, \cite{sq}.
Another instance where the relevance of the $G$--Frobenius algebra is apparent
is in the
tensor product which exists
on the level of $G$--Frobenius algebras and not their invariants.
The philosophy
extends beyond the level of Frobenius algebras to
their deformations, $G$--cohomological field theories as
demonstrated in \cite{JKK}.

The dualization as we described it in \cite{orb} and which we will review
below, does not always provide a $G$--Frobenius algebra. In fact
generally the data of the $\dg$ model with metric does not afford
a $G$--Frobenius algebra structure, although it is expected that there
is a Frobenius structure on the invariants. This leads us to define
the notion of a degenerate $G$--Frobenius algebra below. Here one
adds an additional metric which is equal to the original metric
when restricted to the invariants, but is allowed to be degenerate on
the non--invariant elements and is invariant w.r.t.\
a $G$--graded multiplication. This multiplication
together with the metric descend to a Frobenius algebra on the invariants.

It is this type of structure that arises in the theory of spin
curves \cite{JKV,PV,P} and the construction of cohomological field theories
from certain singularities with fixed Abelian groups $H$ containing
the grading symmetry $J$, which have recently
started to be investigated \cite{FJR}. We conjecture that the resulting
theory is the deformation of the dual of the orbifold of the singularity
with respect to the group $H$.
Although the structures coincide on the
invariant part, on the degenerate part the matching of
non--invariant elements is only almost realized. There are additional
elements which can be explained by interpreting the Milnor rings inside
degenerate $G$--Frobenius algebras, as we discuss in \S \ref{spinparagraph}.


In these geometric settings the $g$--twisted sectors ---which is
another name for the group degree $g$ part of the Frobenius
algebra for $g\neq e$-- which have a degenerate metric, are
related a certain behaviour called of of Ramond type. In the case
of the $A_n$ singularities there is only one such sector and the
entire sector is degenerate. In the cases of $D$ and $E$, the
structure is more complicated and there are invariant elements in
$g$--twisted sectors which have a degenerate metric. The
appearance of these degenerate elements is stunning and maybe a
nuisance from the point of view of spin--curves, but is natural
from the $G$-Frobenius point of view. Moreover regarding our
dualization on the level of $G$--Frobenius algebras as mirror
symmetry, we expect this kind of behavior for the mirror dual
``$A$--model'' of a singularity, the construction of which was
Witten's original motivation for considering the spin--curve
picture \cite{wittenspin}.

A note of caution about nomenclature. One would be inclined
to call the sectors having degenerate pairings in the new metric
Ramond sectors. This might however lead to confusion, since the
term Ramond already has a meaning in the theory of $G$--Frobenius algebras
\cite{orb} and orbifold Landau--Ginzburg theory. Therefore
we will call them sectors of Ramond type and hope to avoid the confusion.

We recall that the Ramond $G$ algebra or state--space
for a $G$--Frobenius algebra is
a cyclic module for the $G$-Frobenius algebra whose $G$--action
is determined by compatibility and
the fact that the generator of the cyclic algebra
is the one dimensional representation of $G$ which is given by the
character $\chi$ which is part of the data of a $G$--Frobenius algebra.
The component of this space of group degree $g$ would be naturally
called the $g$--twisted Ramond sector. The Ramond in this name
stand for the Ramond ground states.
This Ramond space plays a fundamental role in the theory of singularities
as it corresponds as a $\dg$ module to the middle dimensional cohomology
of the Milnor fibers in an orbifold model, while the $G$--Frobenius
algebra corresponds as a $\dg$ module to the orbifold Milnor ring
or universal deformation space. (See the remarks in \S
\ref{singtheory} below). For the untwisted sectors, i.e. the
subalgebras of group degree $e$, this statement
was first proved in \cite{Wa}.

In the sprit of the mirror construction for simple
singularities  one expects that for a given
theory $T$ with a symmetry group $G$ and a subgroup $H\subset G$ of
symmetries $(T/H)^H \simeq (((T/H)/(G/H))^{\vee})^{(G/H)}$ where the subscript
stands for taking the invariants and $\vee$ stands for dualizing.
This type of transformation was used
by \cite{GP} to produce the first mirror pairs. The general statement has to be taken as always
{\it cum grano salis}, but as we show below it is true in many
instances.

Lastly the untwisted sector of an orbifold associated to a singularity
can under certain conditions be related to the folding of an associated
Dynkin diagram. We emphasize that there are foldings and orbifoldings
of diagrams. The $\Zz$ folding of $A_{2n-1}$ yields $B_n$ while
the $\Zz$ orbifolding yields $D_{n+1}$.

In order to understand the operation of folding, we also include
section \ref{folding} in this paper on groups of projective symmetries. This is
a new construction for Frobenius algebras which
we expect to be able to extend to the respective Frobenius manifolds
and to a full theory of $G$--Frobenius algebras.
On the algebra level, we obtain the classical folding results for Coxeter groups
identifying the sub--Frobenius algebra with the Coxeter group of
the folded diagram \cite{strachan}. The relation to singularity theory
and the Milnor fibration is also briefly discussed.
One could hope to extend the folding
to all the diagrams of \cite{zuber} and find the corresponding orbifold theory.


We will work in the setting of $G$--Frobenius algebras over a
field $k$ of characteristic zero (or prime to $|G|$) as it was
established in \cite{orb}. To understand the constructions of
$G$--Frobenius algebra it is important to see that they are
usually performed in four steps. 1. One constructs a $G$--graded
$k$-module $A=\bigoplus A_g$ with a non--degenerate paring between
$A_g$ and $A_{g^{-1}}$ and together with an $A_e$ module structure
on $A$. $A_e$ is usually called the untwisted sector and $A_g$ is
called the $g$--twisted sector. 2. One constructs a  $\dg$ module
structure on $A$ compatible with the $A_e$ module structure. I.e.\
one gives an action of $G$ together with a character $\chi \in
\mathrm{Hom}(G,k^*)$ $\varphi$ s.t.\ $\varphi(g)(A_h)\subset
A_{ghg^{-1}}$ which satisfies the so--called restricted trace
condition and the self--invariance for the twisted sectors. 3.
Lastly one adds a $G$ multiplication to make the $\dg$ module into
a $G$--Frobenius algebra. Comparing $G$--Frobenius algebras on
different levels of this construction compares to the topological
mirror symmetry of dimensions and vector spaces vs.\ that of full
Frobenius manifolds.

These are also the steps of the (re)construction program as explained
in \cite{orb,sq}. Here the data for the first step is usually
provided by the geometric setup. For the second step there are
usually several different choices. This is, however, expected,
since there is the phenomenon of discrete torsion for orbifolds.
As we demonstrated in \cite{disc} for every $G$--Frobenius algebra
there exists a family of $G$--Frobenius algebras indexed by elements
of $\alpha \in Z^2(G,k^*)$ with the same underlying data as mentioned
in step 1 (up to a re--scaling of the metrics pairing the twisted
sectors). In the last step there is an additional
 compatibility condition of the pairing, which
might force one to again re--scale the pairings between the twisted sectors.
For all the conditions, we refer to \cite{orb}. We will however
review the construction for singularities with symmetries below.

The dualization is an involution on triples $(A,j,\chi)$ of a $\dg$ module
$A$, an element $j\in Z(G)$ the center of $G$ and a one--dimensional
representation of $G$, also known as a character. For a special type of graded
$G$--Frobenius algebras $\chi$ is part of the data while $j$ corresponds
to the grading operator. If one includes the other structures of a $G$--Frobenius
algebra, then the operation ceases to be an involution as for instance
the metric will be compatible with the group grading only up to a shift.
To compensate the different behavior of the duals, we introduce the
notion of a degenerate $G$--Frobenius algebra of a group degree $j$ for an
element $j\in Z(G)$.

The paper is organized as follows: In the first section,
we review the construction and basic properties of $G$--Frobenius algebras
and consider special types of graded $G$--Frobenius algebras called Euler
and $G$--Euler.
The second section contains the definition for the dualization for
Euler $\dg$ modules. The third section applies the first two sections
to the $G$--Frobenius algebras resulting from quasi--homogenous polynomials
in general. The fourth section contains  explicit calculations
for a large list of examples. From these examples we obtain
the theorem
about the mirror--self duality of the simple singularities, i.e.\
those of ADE type and the Pham singularities for
coprime powers. The examples also provide mirror pairs for the simple boundary
singularities $B_n$ and $F_4$ and produces $G_2$ as the untwisted sector
of a $D_4$ orbifold.
In the last section, we connect our calculations to spin--curves,
classical results in the theory of singularities and
foldings of Dynkin diagrams.

\section*{Acknowledgements}
We would like to thank the organizers of the AMS special session
``Gromov-Witten Theory
of Spin Curves and Orbifolds''
Tyler Jarvis, Takashi Kimura and Arkady Vaintrob for doing such
an excellent job in providing a common platform for all the
exciting research revolving around spin curves and orbifolds.

It is also a pleasure to thank Boris Dubrovin, Claus Hertling,
 Tyler Jarvis,
Takashi Kimura, Yongbin Ruan, Ian Strachan, Kyoji Saito,
and  Jean--Bernard Zuber for discussions
which were relevant for various stages
of the presented research.

\section{Graded $G$-Frobenius algebras}
\subsection{$G$--Frobenius algebras}
We would like to recall the definition of a $G$--Frobenius algebra
of \cite{orb}. Although it has now appeared in many places we think
it convenient for the reader to display it here once more.

\begin{df}
 A $G$--Frobenius algebra (FA) over
a field  $K$ of characteristic 0 is given by the data
$<G,A,\circ,1,\eta,\varphi,\chi>$, where

\begin{tabular}{ll}
$G$&finite group\\
$A$&finite dim $G$-graded $K$--vector space \\
&$A=\oplus_{g \in G}A_{g}$\\
&$A_{e}$ is called the untwisted sector and \\
&the $A_{g}$ for $g \neq
e$ are called the twisted sectors.\\
$\circ$&a multiplication on $A$ which respects the grading:\\
&$\circ:A_g \otimes A_h \rightarrow A_{gh}$\\
$1$&a fixed element in $A_{e}$--the unit\\
$\eta$&non-degenerate bilinear form\\
&which respects grading i.e. $g|_{A_{g}\otimes A_{h}}=0$ unless
$gh=e$.\\
\end{tabular}

\begin{tabular}{ll}
$\varphi$&an action  of $G$ on $A$
(which will be  by algebra automorphisms), \\
&$\varphi\in \mathrm{Hom}(G,\mathrm{Aut}(A))$, s.t.\
$\varphi_{g}(A_{h})\subset A_{ghg^{-1}}$\\
$\chi$&a character $\chi \in \mathrm {Hom}(G,K^{*})$ \\

\end{tabular}

\vskip 0.3cm

\noindent Satisfying the following axioms:

\noindent{\sc Notation:} We use a subscript on an element of $A$
to signify that it has homogeneous group degree  --e.g.\ $a_g$
means $a_g \in A_g$-- and we write $\varphi_{g}:= \varphi(g)$ and
$\chi_{g}:= \chi(g)$.

\begin{itemize}

\item[a)] { Associativity}

$(a_{g}\circ a_{h}) \circ a_{k} =a_{g}\circ (a_{h} \circ a_{k})$

\item[b)] { Twisted commutativity}

$a_{g}\circ a_{h} = \varphi_{g}(a_{h})\circ a_{g}$

\item[c)] { $G$ Invariant Unit}:

$1 \circ a_{g} = a_{g}\circ 1 = a_g$

and

$\varphi_g(1)=1$

\item[d)] { Invariance of the metric}:

$\eta(a_{g},a_{h}\circ a_{k}) = \eta(a_{g}\circ a_{h},a_{k})$

\item[i)] { Projective self--invariance of the twisted sectors}

$\varphi_{g}|A_{g}=\chi_{g}^{-1}id$

\item[ii)] { $G$--Invariance of the multiplication}

$\varphi_{k}(a_{g}\circ a_{h}) = \varphi_{k}(a_{g})\circ
\varphi_{k}(a_{h})$

\item[iii)]

{ Projective $G$--invariance of the metric}

$\varphi_{g}^{*}(\eta) = \chi_{g}^{-2}\eta$

\item[iv)] { Projective trace axiom}

$\forall c \in A_{[g,h]}$ and $l_c$ left multiplication by $c$:

$\chi_{h}\mathrm {Tr} (l_c  \varphi_{h}|_{A_{g}})=
\chi_{g^{-1}}\mathrm  {Tr}(  \varphi_{g^{-1}} l_c|_{A_{h}})$
\end{itemize}
\end{df}

We sometimes denote by $\rho\in A_e$ the element dual to $\eps\in
A_e^*$ and Poincar\'e dual to $1\in A_e$.

\begin{rm}
For the examples in \S \ref{examples} it is essential that we
consider $G$--Frobenius algebras with non--trivial characters.
\end{rm}

\begin{rk}
Instead of using a left action of $G$ on $A$ one can also
use  a right action as for instance is done in
e.g.\ \cite{JKK}).Since
if $\varphi$ is a left action $\rho(g):=\varphi(g^{-1})$ is a
right action, it does not matter which choice is made.
\end{rk}

\begin{rk}
\label{dgmoduleremark}
Another way to characterize a the $G$--grading and $G$--action
it to say that it is a $\dg$ module.
This statement is equivalent to saying that $A$ is $G$--graded
and the $G$--action is such that $(*)\varphi(g)A_h
\subset A_{ghg^{-1}}$ or $\rho(g)A_h \subset A_{g^{-1}hg}$,
cf.\ e.g.\ \cite{disc}.  We use the
nomenclature of $D(k[G])$ module, rather than $G$--graded $G$--module
since it includes the condition (*).

The compatibilities of the multiplication with the grading and the
$G$--action can also be rephrased as $A$ is a $\dg$ module algebra.
\end{rk}

\subsection{Restriction}
 The operation of
restriction a $G$-Frobenius algebra to a $H$ Frobenius algebra for
a subgroup $H \subset G$ is discussed in \cite{orb} and is given
by  $res(A)^G_H:=\bigoplus_{h\in H} A_ h$ and restricting all
structures.

By forgetting or omitting the multiplicative structure
and considering just the
action of the subgroup $H$ we obtain the restriction
from a $\dg$ to an  $D(k[H])$ module.

\subsection{Super-grading}

We also need to enlarge the framework by considering super--algebras rather
than algebras. This will introduce the standard signs.

\begin{df} A G-twisted Frobenius super--algebra over
a field  $K$ of characteristic 0 is
$<G,A,\circ,1,\eta,\varphi,\chi>$, where \vskip 0.3cm

\begin{tabular}{ll}
$G$&finite group\\
$A$&finite dimensional  $\Zz \times G$-graded $K$--vector space \\
&$A= A_0 \oplus A_1= \oplus_{g \in G}(A_{g,0}\oplus A_{g,1}) = \oplus_{g \in G}A_{g}$\\
&$A_{e}$ is called the untwisted sector and is even. \\
&The $A_{g}$ for $g \neq  e$ are called the twisted sectors.\\
$\circ$&a multiplication on $A$ which respects both gradings:\\
&$\circ:A_{g,i} \otimes A_{h,j} \rightarrow A_{gh,i+j}$\\
$1$&a fixed element in $A_{e}$--the unit\\
$\eta$&non-degenerate even bilinear form\\
&which respects grading i.e. $g|_{A_{g}\otimes A_{h}}=0$ unless
$gh=e$.\\
$\varphi$&an action by even algebra automorphisms  of $G$ on $A$, \\
&$\varphi\in \mathrm{Hom}_{K-alg}(G,A)$, s.t.\
$\varphi_{g}(A_{h})\subset A_{ghg^{-1}}$\\
$\chi$&a character $\chi \in \mathrm {Hom}(G,k^{*})$
\end{tabular}

\noindent satisfy the  axioms a)--d) and i)--iii) of a
$G$--Frobenius algebra with the following alteration:

\begin{itemize}

\item[b$^{\sigma}$)] { Twisted super--commutativity}

$a_{g}\circ a_{h} = (-1)^{\tilde a_g\tilde a_h}
\varphi_{g}(a_{h})\circ a_{g}$

\item[iv$^{\sigma}$)] { Projective super--trace axiom}

$\forall c \in A_{[g,h]}$ and $l_c$ left multiplication by $c$:

$\chi_{h}\mathrm {STr} (l_c  \varphi_{h}|_{A_{g}})=
\chi_{g^{-1}}\mathrm  {STr}(  \varphi_{g^{-1}} l_c|_{A_{h}})$
\end{itemize}
where $\mathrm{STr}$ is the super--trace.
\end{df}

\subsection{Graded $G$--Frobenius algebras}

\begin{df}
We call a (super) $G$ Frobenius algebra $A$ graded by an additive
group $I$ if it is graded as a (super) algebra by $I$ and the metric
is homogenous of a fixed degree $d$, i.e.\ for homogenous $a,b$,
$\eta(a,b)=0$ unless $deg(a)+deg(b)=d$, where we denote
the $I$ degree of a homogenous element $a\in A$  by $\deg(a)$. If
$I=\mathbb{Q}$, we simply call $A$ graded.
We also call $d$ the degree of the Frobenius algebra.
\end{df}

The degree of the Frobenius algebra is the degree of the element $\rho$.


\subsection{The grading operator}
Given a graded $G$--Frobenius algebra $A$, we define the grading
operator $\Q$ to be given by

\begin{equation}
\Q(a) := \deg(a)a\quad \text{if $a$ is homogeneous}
\end{equation}

Sometimes this type of operator is also called $E$.

In the case that $A$ is graded and $k= \mathbb{C}$ or $k$ is of characteristic
$0$ and an embedding
of $\bar k \subset \mathbb{C}$ has been fixed
we furthermore define the operator

\begin{equation}
J:= \exp(2\pi i \Q)
\end{equation}

\begin{df}
We call a graded $\dg$--module $A=\bigoplus_{g \in G} A_g$
Euler if the operator $J|_{A_e}$
is described by the action of a central element $j$ of the group $G$ on $A_e$.
I.e. there exists a $j \in Z(G)$, the center of $G$,
$\varphi(h^{-1}j)|_{A_h}= J|_{A_h}$.

We call a graded $\dg$--module $G$--Euler if
 there exists a $j \in Z(G)$, s.t.\
$\varphi(h^{-1}j)|_{A_h}= J|_{A_h}$.

We call a graded $k[D(H)]$ $A$ quasi--Euler (or quasi--$G$--Euler)
if there is a group $G$, s.t. $H$ is
a subgroups of $G$ ($G\supset H$) and there exists an Euler (or $G$--Euler) $\dg$
module $B$ s.t.\ the restriction of the $\dg$ module $B$
to its $k[D(H)]$--sub module  $res_H(B)$ is $A$.

An Eulerization of a quasi--Euler $D[k(H)]$ module is a fixed choice of $\dg$ module
$B$ as above.

A $G$--Frobenius algebra is called Euler, $G$--Euler, quasi--Euler or quasi--$G$--Euler
if its underlying $\dg$ module is Euler, $G$--Euler, quasi--Euler or quasi--$G$--Euler,
respectively.
\end{df}

\subsection{Bi--Graded $G$--Frobenius algebras}

\begin{df}
We call a (super) $G$ Frobenius algebra $A$ bi--graded by an
additive group $I$ if it is bi--graded as a (super) algebra by
$I$. If $I=\mathbb{Q}$, we simply call $A$ bi--graded.
\end{df}

\subsubsection{Notation} We will usually denote the two grading operators by
$\Q$ and $\bar \Q$. Given a bi--homogenous element $a$ we will
denote its degree w.r.t.\ $\Q$ by $q(a)=\Q(a)$ and its degree
w.r.t.\ $\bar \Q$ by $\bar q(a)=\bar \Q(a)$. We will also use the
notation $(q(a),\bar q(a))$ to denote the bi--degree.

\begin{df}
Fix a graded $G$--Frobenius algebra $A$ with grading operator
$\mathcal{Q}$.

We define its $(c,c)$ realization
$A^{(c,c)}$ to be given by the $G$--Frobenius algebra $A$ together with the
bi--grading $(\mathcal {Q},\mathcal {Q})$, i.e.\
$\bar{\mathcal{Q}}=\mathcal{Q}$.

We define the $(a,c)$ realization of $A$ denoted by $A^{(a,c)}$
to be given by the $G$--Frobenius algebra $A$ together with the
bi--grading $(\mathcal {Q},-\mathcal {Q})$,
i.e.\ $\bar{\mathcal{Q}}=-\mathcal{Q}$.

\end{df}


\begin{rk}
The terminology stems from the representation theory of the $N=2$ super--conformal
algebra, as explained in the introduction.
\end{rk}

\subsection{Constructing $G$--Frobenius algebras}
\label{conststeps}

When constructing $G$--Frobenius algebras form geometric or
algebraic data usually the different structures are introduced one
after the other.  A good example of this procedure is given by the
construction of $G$--Frobenius algebras from isolated
singularities with symmetries reviewed below
\ref{singconstruction}. Also some operations like the duality
discussed below  are given on a certain level of structure. The
usual order in which the structures are introduced is as follows.

\begin{enumerate}

\item {\bf The $G$--graded $k$--module.} Usually the first
structure to be given for any $G$--Frobenius algebra is its {\em
additive} structure $A:=\bigoplus_{g\in G} A_g$.

On this level it is also usual to introduce the non--degenerate
pairing $\eta$ which pairs $A_g$ with $A_{g^{-1}}$.

\item {\bf The $G$--graded $G$--module or $D(k[G])$ module
structure.} The next property which is usually introduced is a
{\em $G$--action} on $A$, usually denoted by $\varphi$ for a left
action (cf.\ e.g.\ \cite{orb}) which makes $A$ into a $\dg$--module
cf.\ \ref{dgmoduleremark}.

Further data and conditions:
\begin{enumerate}
\item Along with the $G$--action the function $\chi:G\rightarrow
k^*$ is fixed
since $\chi$ can be derived from the $G$--action via the condition
of projective self--invariance (axiom ii$^{\sigma}$).

\item From the projective $G$--invariance of the metric (axiom
iii), it follows that the function $\chi^2$ has to be a character,
i.e. a one--dimensional representation.

\end{enumerate}

\item {\bf The  $D(k[G]),A_e$--bi--module} Usually the
untwisted sector $A_e$ is naturally a Frobenius algebra. The next
step in constructing a $G$--Frobenius algebra is then the {\em
$A_e$--module} structure for each $A_g$: $A_e\otimes A_g
\rightarrow A_g$, which
 will be a part of the algebra  multiplication.
 These operations turn $A$ into an $A_e$ module. They are usually
 already present in the geometry by functoriality \cite{sq}.
 The $A_e$ module structure
 has be compatible with the $G$--action so: $\varphi(g)(a_e b_h)=
 \varphi(g)(a_g)\varphi_g(b_h)$.

The $A_e$ module structure leads to a second compatibility
condition of the $G$--action with the pairing which is given by
restriction of the trace axiom to the case $g=e ,c=1 \in A_e$.
This condition effectively relates the dimension of the various
twisted sectors $A_g$ to the character $\chi$ and $G$--action on
the identity sector.
$$
\chi_h \mathrm{STr} \varphi(h)|_{A_e}= \mathrm{STr} id|_{A_h}=
\mathrm{sdim}(A_h)
$$
Also the trace axiom put constraints on the possible $G$--actions.
The constrains can be quite effective, but they define
the action at most up to discrete torsion \cite{disc}.

\item{\bf The $G$--Frobenius algebra.} The last step is to
introduce the {\it stringy multiplications}:$A_g \otimes A_h
\rightarrow A_{gh}$ , ie. the algebra structure. This structure
has to be compatible with the $G$--action and the metric.

\end{enumerate}

\subsection{The metric and the grading}
\label{metricandgrading} When constructing a
$G$--Frobenius algebra in the above fashion, the metric and the
grading can either be introduced at the end, but usually, there is
a natural choice in each step, which may be modified in the next
step.

\subsubsection{The metric}
The metric, i.e. non--degenerate even symmetric pairing, is
usually introduced in step (1) and may be re--scaled in step (4)
by a factor to ensure the compatibility of the metric with the
multiplication (invariance of the metric axiom d)).

\subsubsection{The grading}
In step (1) first there is usually a grading $\Q^{(1)}$ inherent in the
definition of each $A_g$ when introducing the metric which is
usually also inherent in the construction. Each of the pairings
$A_g \otimes A_{g^{-1}}\rightarrow k$  is
usually homogenous of some fixed degree $d_g$ with
$d_g=d_{g^{-1}}$.

The first alteration of the grading is a shift of the grading for
each $A_g$ by $\frac{1}{2}s^+(g):= \frac{1}{2}(d_e -d_g)$, i.e.\
the new grading for an element $a_g \in A_g$ is  $\Q^{(2)}(a_g)= \Q^{(1)}(a_g)
+ \frac{1}{2}s^+(g)$. This makes the metric homogenous of degree
$d_e$ on all of $A$. Notice that $s^+(g)=s^+(g^{-1})$.

The second alteration appears in step (2). For physically inspired
reasons, one often  makes an additional shift $\frac{1}{2}
s^{-1}(g)$ depending on $g$  which has to preserve the
homogeneity of the metric. The second shift satisfies
$s^{-1}(g)=- s^{-1}(g^{-1})$.

The final grading for an element $a_g
\in A_g$ is
$$\mathcal {Q}(a_g)= \Q^{(2)}(a_g) + \frac{1}{2}s^-(g)=  \Q^{(1)}(a_g) +
\frac{1}{2}(s^+(g)+s^-(g)).$$

If the $G$--action of step (2) is
induced by a linear $G$ action there is a standard choice for this
shift, given by

\begin{df}
The standard grading shift for a $G$--Frobenius algebra with a
choice of linear representation $\rho: G \rightarrow GL_n(k)$ is
given by
\begin{equation}
\label{sdef} s_g := \frac{1}{2}(s_g^+ + s_g^-)
\end{equation}
with
\begin{equation}
\label{s+def} s_{g}^+:= d-d_{g}
\end{equation}
 and
\begin{eqnarray}
\label{s-def}
 s_{g}^- := \frac{1}{2\pi i}\mathrm{tr}
(\log(g))-\mathrm{tr}(\log(g^{-1}))&:=& \frac{1}{2\pi i}(\sum_i
\l_i(g)-\sum_i \l_i(g^{-1}))\nn\\
 &=&\sum_{i: \l_i \neq 0}
(\frac{1}{2\pi i}2\l_i(g)-1)
\end{eqnarray}
where the $\l_i(g)$ are the logarithms of the eigenvalues of
$\rho(g)$ using the arguments in $[0,2\pi)$.

This means that if $\rho(g)= diag(exp(2\pi i \nu_1),
\dots,exp(2\pi i \nu_n))$ with $0\leq \nu_i< 1$ then $\l_i = 2\pi
i \nu_i$.
\end{df}

\begin{rk}
Notice that if $\rho(g)= diag(exp(2\pi i \nu_1(g)), \dots,
exp(2\pi i \nu_n(g)))$ with $0\leq \nu(g)_i< 1$ then
$$
\nu_i(g^{-1}) = \begin{cases} 0 &\text{if } \nu_i(g) =0 \\
1-\nu_i(g) & else \end{cases}.
$$
$$
\nu_i(gh) = \nu_i(g) + \nu_i(h) - \Theta(1- (\nu_i(g) + \nu_i(h))
$$
where is the step function
$$\Theta(x) =\begin{cases} 1
&\text{if } x \geq 0\\0 &\text{if } x<0\end{cases}$$
\end{rk}

\begin{rk}
In the case of orbifold cohomology \cite{CR}, one starts with an
action of $G$ on the manifold $M$ and induces an action on the tangent
space $M$ which defines the shift $s^-$ via (\ref{s-def}) and the
shift $s^+_g$ is defined by $d_g:= \dim (Fix(g)\subset $M$)$. For
general orbifolds this reasoning is understood locally \cite{CR}.
For global orbifolds the expressions can however be understood
globally.

If $\rho(g)= diag(\exp(2\pi i \l_1), \dots,\exp(2\pi i \l_n)$
,then $d_g=\sum_{i:\l_i=0} 1$ and $d-d_g = \sum_{i:\l_i \neq 0}
1$, so $s_g =\sum_{i: \l_i \neq 0}\frac{1}{2\pi i}\l_i(g)= \sum_i
q_i$ yielding agreement with the definition (\ref{sdef}) above and
the one of \cite{CR} and \cite{zuber} in that particular case.
\end{rk}

Notice that for the last expression of equation (\ref{s-def}), we
can use the branch of the logarithm obtained by cutting along
$[0,\infty)$.
\subsubsection{The super--grading} As for the grading, usually
each $A_g$ comes with an intrinsic super--grading. In step (1) one
usually allows the freedom to shift the super--grading by
$\mathbb{Z}/2\mathbb{Z}$ values function. The restrictions on this
function come from the existence of an even non--degenerate
quasi--homogenous metric and in step  (3) from the trace axiom. In
step (4) the condition that $\chi$ is a character translates via
the trace condition and the condition that $\chi^2$ is a character
from step (2) into a condition on the super--grading.

\subsubsection{Bi--grading}

\begin{df}
Set $\bar s_g:=\frac{1}{2}(s^+_g-s^-_g)$. Since $s^+_g=s^+_{g^{-1}}$ and
$s^-_g=-s^-_{g^{-1}}$ if follows that $\bar s_g=s_{g^{-1}}$.
 We define that bi--grading $(\cQ,\cbQ)$ by
$$
 \cQ(a_g):= \Q(a_g)+ s_g \quad \cbQ(a_g):= \Q(a_g)+ \bar s_g \quad
 \text{ for } a_g \in A_g=\M_{f|_{Fix(g)}}
$$
\end{df}

\begin{rk}
As mentioned in the introduction, the bi--grading has its origin
in the interpretation for the
algebra as the $(c,c)$ ring for an orbifold model \cite{IV}.
\end{rk}

\section{A mirror type transformation}

{\bf Assumptions:} In this section for simplicity, we fix a $k=
\mathbb {C}$. (If $k$ is a field of characteristic zero, we could
fix an embedding $\bar{k} \hookrightarrow \mathbb{C}$.)

In the following, we will construct an involution for the
triples $\langle A,j,\chi\rangle$ of
$D(k[G])$--modules $A$, elements $j$ of the center of $G$ and
characters $\chi \in \mathrm{Hom}(G,k^*)$.

In the case of an Euler $\dg$ module, we take the element $j$ to be
the element defined by the Euler property.

We also extend the operation to include a non--degenerate pairing
and a bi--grading.

This involution induces via restriction a dualization on
quasi--Euler  $\dh$ modules (without pairing)
with fixed Eulerization.

In the case $A$ is a Euler Frobenius algebra or
a quasi--Euler Frobenius algebra with a fixed Eulerization,
we let $\Q$ be the grading operator and $j \in G$, s.t. $\rho(j)=
exp(2\pi \Q)=J$. In this case the data $(A,j,\chi)$ is fixed
by the $G$--Frobenius algebra and the element $j$ yielding the grading.

Our operation conjecturally
acts as a mirror transformation on the underlying Euler $G$--Frobenius algebras
in the sense of \cred, see \S \ref{credo}.

The  additional bi--grading, is conjecturally
compatible with interchange of the $(c,c)$-type and $(a,c)$-type for
Landau--Ginzburg theories w.r.t.\ mirror symmetry.

In fact, we will prove that the \credo is correct in the case of the
simple singularities $A_n,D_n,E_6,E_7,E_8$ and
yields mirror pairs for the simple boundary
singularities $B_n,F_4$.

\begin{rk} The definition of the dual comes from physics
\cite{IV,V}, where the dual $D(k[G])$--module is obtained by using
an endofunctor in the category of representations of the $N=2$
super--conformal algebra which translates in our case to an
isomorphism of $D(k[G])$--modules. This endofunctor is generally
known as spectral flow and has a particular realization discussed
below in the case of $G$--Frobenius algebras.
\end{rk}

\subsection{The $G$--graded $k$--module structure}

\begin{df}
 Given a $G$--graded $k$--module $A$ and an element $j\in Z(G)$
 we define the dual $\check A$ to be the $G$--graded $k$ module:
 \begin{equation}
 \check A_g := A_{gj^{-1}}, \quad \check A := \bigoplus_{g \in G}
 \check A_g
 \end{equation}
\end{df}

\begin{rk}
The above formula states that as $k$ modules  $A$ and $\check A$
are isomorphic. It is only their $G$--grading which has changed.
We denote the isomorphism by $M: A \rightarrow \check A$, with
$M(A_g)=\check A_{gj}$.
\end{rk}

\subsection{The metric}
With the help of the map $M^{-1}$, we can pull back a given metric
$\eta$ from $A$ to $\check A$. We set

$$
\check \eta = (M^{-1})^{*}\eta \qquad \check\eta(\check a,\check
b) :=\eta(M^{-1}(\check a),M^{-1}(\check b))
$$

\begin{rk}
\label{metrk} Notice if $\eta$ is homogeneous with respect to the
group degree, i.e. pairs $A_g$ with $A_{g^{-1}}$,
then $\check \eta$pairs $\check A_g$ with $\check
A_{g^{-1}j^{2}}$ and thus $\check \eta$ is not group degree
homogeneous, but of group degree $j^{2}$ as a tensor in
$\check A^*\otimes \check A^*$.
\end{rk}

\begin{rk}
The metric $\check \eta$ is $G$--invariant and hence descends
to the $G$--invariants.
\end{rk}

\subsection{The $G$--action or the $D(G[k])$--module structure}

Given a triple $\langle A,j,\chi\rangle$ of a
$D(k[G])$--module $A$, an element $j$ of the center of $G$ and
a character $\chi \in \mathrm{Hom}(G,k^*)$, we
define $\bar \varphi:= \varphi \otimes_k \chi$.
This is an action of $G$ on the
$k$--module $A\otimes_k k \simeq A$ and thus on the $k$--module $\check A$.

We define the $G$-action $\check \varphi$ on $\check A$ to be
the induced action of the action on $A$ by $\bar \varphi$. That is for
$a\in \check A$

\begin{equation}
 \check {\varphi}(h)(a) = \chi(h) M (\varphi(h)(M^{-1}(a)))
\end{equation}

\begin{rk}
We see that under this action $\check\varphi (h)(\check A_g)
\subset \check A_{hgj^{-1}h^{-1}j}=A_{hgh^{-1}}$, since we made
the assumption that $j\in Z(G)$. Thus we obtain a $G$--action,
which makes $\check A$ into a $D(k[G]])$--module.
\end{rk}

\subsection{The bi--grading of the dual.}
If $A$ was initially graded by the operator $Q^{(1)}$ and
or simply
$Q(a_g)=Q^{(1)}(a_g) +s_g$ then.

Set $\check s_g:=s_{gj^{-1}}-d$ and $\check {\bar s}_g:=\bar
s_{gj^{-1}}$, where we recall that $d$ is the
degree of the $G$--Frobenius algebra.
We define a bi--grading on $\check A$ by
\begin{equation}
\label{dualbigrading} \check {\mathcal{Q}}(\check a) = \Q^{(1)}(a)+\check
s_g \quad \bar{\check{\mathcal{Q}}} := \Q^{(1)}(a) +\bar{\check s}_g
\quad \text{ for } \check a_g \in \check A_g
\end{equation}

\begin{rk}
For an Euler $G$--Frobenius algebra $A=<G,A,\circ,1,\eta,\varphi,\chi,j>$
naturally gives rise to a triple
$<A,j,\chi>$ and thus obtain a dual $\dg$ module with a non--degenerate pairing
and a bi--grading.
\end{rk}

\begin{rk}
The motivation for the dual bi--grading again comes from the physical
interpretation of $G\M_f$ as an orbifold Landau--Ginzburg model
and the dualization being implemented by the spectral flow
operator $\mathcal{U}_{(1,0)}$ \cite{IV} which has
the natural charge $(d=\hat c=\frac{c}{3},0)$.
\end{rk}

\subsubsection{The involution}

\begin{df}
We define the dual of a triple $\langle A,j,\chi\rangle$ of
$D(k[G])$--modules $A$, elements $j$ of the center of $G$ and
characters $\chi \in \mathrm{Hom}(G,k^*)$ to be the triple
$\langle \check A, j^{-1}, \chi^{-1} \rangle$.
\end{df}

\begin{rk}
Notice that the inclusion of the data $j$ and $\chi$ turns the operation
on the $\dg$--module into an involution.
\end{rk}

\subsubsection{The dual of a quasi--Euler $\dh$--module with given Eulerization}

\begin{df}
We define the dual of a quasi--Euler $\dh$--module $A$ (or $H$--Frobenius algebra)
with given Eulerization $B$
to be the restriction of $\check B$ to $H$. $\check A:= res_H(\check B)$
\end{df}

\begin{rk}
Notice that if $j\notin H$ then we cannot pull back the metric, since
if $h\in H$, $h^{-1}j^{2}$ need not be in $H$. If $\forall h \in H: h^{-1}j^{2} \in H$,
then we can also pull back the metric.
\end{rk}

\subsection{A dual $G$--Frobenius algebra?}
We would like to remark that the dualizing process is only a
process of dualizing for $D(k[g])$ modules with metric.

One thing
to prevent the resulting structure from being a $G$--Frobenius
algebra is that the metric is not $G$--graded anymore as remarked
above in Remark \ref{metrk}. Also the projective self--invariance
might not hold. However, there might be, in some cases unique,
choices of $G$--graded multiplication compatible with the
$G$--action.

Or what is actually expected by physics, that there is a Frobenius
algebra structure on the $G$--invariants of $\check A$ with the
given metric. It is important to note that physics does not say
there should be an algebra isomorphism and in fact the induced
multiplication $M\circ M^{-1}$ will not be $G$--graded on $\check
A$ unless $j=e$ and the grading and dualization are trivial.

What we can expect is a Frobenius structure on the invariants, plus
a lift of this Frobenius structure to the $G$--graded
equivariant level. This will provide
some additional structure. This motivates the following definition.

\begin{df}
A degenerate $G$--Frobenius algebra  $A$ of degree $j\in Z(G)$ is given by
the data $\langle G,A,\circ,1,\eta,\eta',\varphi,\chi \rangle$
where $<G,A,\circ,1,\eta,\varphi,\chi>$ are the data of a
$G$--Frobenius algebra, and $\eta'$ is a second pairing on
$G$. These data  satisfy the conditions of
a $G$--Frobenius algebra with the following changes and additions:

\begin{itemize}
\item[1)] The non--degenerate paring
$\eta$ and the pairing  $\eta'$ pair $A_{gj^{-1}}$ with $A_{g^{-1}j^{-1}}$.
\item[2)] $\eta|_{A^G} = \eta'_{A^G}$ where $A^G$ are the $G$
invariants of $A$.

\item[d')] Invariance of the metric $\eta'$:

$\eta'(a_{g},a_{h}\circ a_{k}) = \eta'(a_{g}\circ a_{h},a_{k})$

\item[i)$^j$] {Self--invariance of the twisted sectors}

$\varphi_{gj^{-1}}|A_{g}=id$

\item[iii$^j$)]

{$G$--invariance of the metric $\eta$}

$\varphi_{g}^{*}(\eta) = \eta$

\item[iv$^j$)] {$j$ twisted trace axiom}

$\forall c \in A_{[g,h]}$ and $l_c$ left multiplication by $c$:

$\mathrm {Tr} (l_c  \varphi_{hj^{-1}}|_{A_{g}})=
\mathrm  {Tr}(  \varphi_{g^{-1}j} l_c|_{A_{h}})$

\end{itemize}
\end{df}

\begin{conj}
We conjecture that there is a degenerate $G$--Frobenius algebra
of degree $j$ on the dual of a $G$--Euler $G$--Frobenius algebra.
\end{conj}

In the examples we consider, there is a certain uniqueness
in the choice for the multiplication. In order to state this
precisely we need the following two definitions.

\begin{df}
Fix a $\dg$, $A_e$ bi--module $A=\bigoplus A_g$ together with
two metrics $\eta$, $\eta'$ and a character $\chi$
satisfying all the axioms pertaining
to the metrics and the $G$--actions of a degenerate $G$--Frobenius
algebra of degree $j$. We call a degenerate $G$--Frobenius
algebra $A$ of degree $j$ maximally non--degenerate if
$a_g \circ b_h =0$ in $A$ implies that $a_g \circ' b_h=0$ in
any other degenerate $G$--Frobenius
algebra $A'$ with the same underlying $\dg,A_e$ bi--module
$A'=\bigoplus A_g$ together with the
two metrics $\eta$, $\eta'$ and the character $\chi$ .

We call a maximally non--degenerate $G$--Frobenius algebra
projectively unique if it agrees with all
other maximally non--degenerate $G$--Frobenius
structures when projected to $A/k^*$.
\end{df}

\begin{rk}
As demonstrated in \cite{disc} twisting by discrete torsion
exactly realizes the universal (i.e. applicable to all $G$--Frobenius algebras)
projective rescalings. This means vice versa that $A$ and all its twists by
discrete torsion are projectively the same.
\end{rk}

\subsection{\Credo}
\label{credo}

There is a \credo  which is motivated by physics (cf. eg.\ \cite{GP})
or representation theory which states the following
\begin{phil}
Let  $T$ be a $N=2$ theory (which for us at the moment means Frobenius algebra)
and let $H\subset G$ be symmetry groups with $H$ normal in $G$ then
$$
(T/H)^H \simeq (((T/H)/(G/H))^{\vee})^{(G/H)}
$$
for us $T/K$ means a $K$--Frobenius algebra derived from $K$.
\end{phil}

This is too vague to be called a conjecture, since most of the
symbols in the statement have no fixed meaning. We can
however apply it to quasi--homogenous singularities, where up
to a finite amount of data ($\sigma,\eps,\gamma)$ if $G$ is Abelian
see \S \ref{singconstruction} below) the ingredients $T/K$ are fixed.

It turns out that even for different actions of $H$ and $G$ this
\credo holds true.

In order to elucidate the statement, we wish to point out that
there is indeed an action of $G/H$ on the $H$ invariants of a $G$--Frobenius
algebra $A$.
For the action to be defined on the restriction  $res_H(A)$ we need
the addition assumption that $H$ is normal. Therefore the statement
makes sense on the level of $\dg$--modules once the $G$--action is fixed.

\section{Quasi--homogenous singularities with symmetries}
\label{singconstruction}

We again fix $k=\mathbb{C}$.

\begin{df}
 Let $f:{\mathbb C}^n
\rightarrow {\mathbb C}$ be a  function  which has an isolated
singularity at zero. A symmetry of $f$ is an element $S\in
GL(n,\mathbf{C})$, s.t.\ $f(S(\mathbf{z}))=f(\mathbf{z})$. An
isolated singularity with symmetries is a function $f:{\mathbb
C}^n \rightarrow {\mathbb C}$  which has an isolated singularity
at zero together with a finite group $G$ and a representation
$\rho(G): G \rightarrow GL(n,\mathbb{C})$ such that $G$ acts by
symmetries on $f$, i.e.\ $\forall g\in G: g^*(f)(\mathbf{z})=
f(\rho (g)(\mathbf{z}))=f(\mathbf{z})$.

We denote by
$G_{max}\subset GL(n,\mathbb{C})$ the maximal group of symmetries
of $f$.
\end{df}


\begin{dfprop}
For a  function $f{(\bf z})$ with an isolated singularity at zero.
We will denote by $\M_f$ the Milnor or local ring of $f$, which is
given by $\mathcal{O}/J_f$ where $\mathcal{O}$ is the ring of germs
of holomorphic functions at zero and $J_f =(\frac{\partial
f}{\partial z_i})$ is the Jacobian ideal. This ring together with
the Grothendieck residue paring $\eta$ is a graded Frobenius
algebra, see e.g.\ \cite{arnoldbook,maninbook}.
\end{dfprop}

\subsection{The graded Frobenius algebra of a quasi--homogenous with an
isolated singularity at zero}

If the function $f$ is also quasi--homogenous there is a natural
grading operator which assigns to each $z_i$ its degree of
quasi--homogeneity $q_i$.

To define the $q_i$ assume that
$$
f(\lambda^{Q_1}z_1, \dots, \lambda^{Q_n}z_n) = \lambda^N
f(z_1,\dots, z_n)
$$
with $Q_i,N \in \mathbb{N}$. Then we set $q_i=\frac{Q_i}{N}$ and
define $\deg(z_i):=q_i$ which yields a map $\M_f \rightarrow
\mathbb{Q}$.

The metric for the resulting
Frobenius algebra is given by the element which is dual to the
identity and this element is represented by $H:=Hess(f)$ the Hessian of $f$.
For a quasi--homogenous singularity,  the degree of the Hessian is
the degree of the form $\eta$ and is denoted by $d$. By the
general theory \cite{arnoldbook} a formula for $d$ is given by
$$
d=\sum_i (1-2q_i)
$$
and the dimension or the Milnor number of the local algebra is
$$
\mu:=\dim(\M_f)=\prod_i (\frac{1}{q_i}-1)
$$

\subsubsection{Examples} In the following examples, we took the
liberty to re--scale the Grothendieck residue form, which amounts to
adding an overall factor to the function $f$.

\begin{enumerate}
\item The $A_n$ series: $f(z,w)=z^{n+1}$
$$\M_f =
\mathbb{C}[z]/(z^{n-1})=\langle 1,z,z^2,\dots,z^{n-1}\rangle  \quad
\eta(z^i,z^j)= \delta_{i+j,n-1}$$ $q=\frac{1}{n+1},
\mu=n,d=1-\frac{2}{n+1}=\frac{n-1}{n+1}$.
\item The $D_n$ series: The $D_{n+1}$, $n\geq 3$ singularity is
given by the function $f(x,y,w)=\frac{1}{n}x^n+ xy^2$
$\M_f \simeq
\mathbb{C}[x,y]/(x^n,xy)\simeq\langle 1,x,x^2,\dots, x^{n-2},y\rangle\\
\eta(x^i,x^j)= \delta_{i+j,n-1}, \eta(y,y)=1,
\eta(x^i,y)=0 \\
q_x=q_1=\frac{1}{n},
q_y=q_2=\frac{n-1}{2n},\mu=n,d=\frac{n-1}{n}$

\item The $E_7$ singularity:$\frac{1}{3}x^3+xy^3$

$\M_f =
\mathbb{C}[x,y]/(x^3,xy^2)=\langle 1,x,x^2,y,y^2,xy,x^2y,
x^{n-2},y\rangle \\
 \eta(x^iy^j,x^ky_l)=
\delta_{i+k,2}\delta_{j+l,1}\\
q_x=q_1=\frac{1}{3}
q_y=q_2=\frac{2}{9}, \mu=6,d=\frac{8}{9}$.
\end{enumerate}

\subsubsection{Products}
 For two functions $f$ and $g$
with an isolated singularity at zero, as shown in \cite{tenfrob,
maninbook} $\M_{f+g}=\M_{f} \otimes \M_g$ even on the level of
Frobenius manifolds.

\subsubsection{Stabilization}
Notice that adding squares , an operations known as
stabilization,
to a function with an isolated singularity ($f\mapsto f+w^2$)
leaves the Milnor ring
invariant. This fact which is well known in singularity theory
(see e.g.\ \cite{arnoldbook}), can also be seen as follows from the
point of view of Frobenius algebras.

Since the Frobenius algebra of the singularity $f(w)=w^2$
 $\M_{w^2}=A_1=k$ is the unit in the monoidal category
 of Frobenius algebras \cite{tenfrob} we also find
 that $\M_{f+w^2}\simeq \M_f\otimes A_1\simeq \M_f$.

In the following, all the definitions, calculations and operations
are invariant under stabilization.

\begin{df}
We define $\M_0:=\M_{w^2}=A^1$.
\end{df}

\begin{rk}
All the following definitions and constructions are invariant under stabilizations,
if one extends the group action by the usual embedding of $GL(n,\mathbb{C})$
to $GL(n+1,\mathbb{C})$.
\end{rk}

\subsection{The $G$--Frobenius algebra for a singularity with symmetry group $G$}
\label{QMf}
We would like to recall from \cite{conf,orb,sq,hilb} that for the
data $(f,G,\rho)$ as above there are several natural
$G$--Frobenius algebras, whose underlying $k$--module structure
and bi--grading are all the same, but whose $D(k[G])$ module
structures are in one---to---one correspondence with twists by discrete
torsion and whose $G$--Frobenius structures depend on the choice
of a graded compatible co--cycle for the quantum multiplication.
We will review the construction below following the steps of \S
\ref{conststeps}.

\subsubsection{The $G$--graded $k$--module structure}
First we show that for the data $(f,G,\rho)$
as above there is a natural associated
$G$--graded $\M_f$ module.

Let $\mathrm{Fix}_g:=$ the fixed point set of $g$ in
$\mathbb{C}^n$, in other words the eigenspace to the eigenvalue
$1$ of $\rho(g)$. Set $f_g:= f|_{\mathrm{Fix}_g}$.

We define
$$
A_g := \M_{f_g} \quad G\M_f := \bigoplus_{g\in G} \M_{f_g}
$$

\begin{rk} We would like to emphasize the following observations:
\begin{itemize}
\item[1)] Notice $A_e= \M_f$.

\item[2)] Each of the $A_g=\M_{f_g}$ is as a local ring of a
quasi--homogenous function with an isolated singularity at zero is a
Frobenius algebra. We denote the metric for the Frobenius algebra
$A_g$ by $\eta_g$, its unit by $1_g$ and its degree by $d_g$.
and its grading operator by $Q_g$.

These sum of the grading operators $Q_g$ defines a grading operator
$Q$ on $A$.

\item[3)] We furthermore can use the ring structure of the
individual $\M_{f_g}$ to define natural $\M_f$ module structure
by inclusion of function germs. This $A_e$--module structure is compatible
with the grading by $Q$.

In the following examples the multiplication is given by
\footnote{
The precise condition is as follows. Let's suppose $g$ is diagonal
in the variables $z_i$ and let $z_i:i \in I_1$ be a basis of $Fix(g)$
and $z_i:i\in I_2$ a basis of the complement of the fixed point set.
Set $Z=(z_i:i \in I_2)$ let $J_g = J_{f_g}$
and $J_{I_1}=(f_{z_i}:i\in I_1)$ be the respective ideals in
$\mathcal{O}$. Then the condition is
that $J_{I_1}+Z = J_g+Z$. In the other case the map still exists by
functoriality, but is a little more complicated.}:

$$
\M_f \times \M_{f_g}\rightarrow \M_{f_g}: (a,b) \mapsto
a|_{\mathrm{Fix}_g} b
$$
\end{itemize}
\end{rk}

\begin{rk}
All the $A_g$ are cyclic $A_e$ modules.
In the terminology of \cite{orb} $G\M_f$ is a special
$G$--Frobenius algebra. Notice that  the  unit  $1_g$ is
a cyclic generator for the $A_e$ module  $A_g$.
\end{rk}

\subsubsection{The grading}

The initial grading operator $Q$ from above plays the role of the operator
$Q^{(1)}$ of \S \ref{metricandgrading}.
The actual grading $\Q$ is determined by the degrees of the cyclic generators
$1_g$.

\begin{df} We define the grading operator $\mathcal{Q}$ on $G\M_f$
by
$$
\mathcal{Q}(a_g)= \Q(a)+s_g \text{for } a_g=a1_g
$$
with
$$
s_{g}= \frac{1}{2}(s_g^+ + s_g^-)= \frac{1}{2}(d-d_g) + \sum_{i:\nu_i\neq 0}
(\frac{1}{2\pi i}\l_i(g)-\frac{1}{2}
)$$

where the $\l_i(g)$ are the logarithms of the eigenvalues of $g$
using the branch with arguments in $[0,2\pi)$ i.e.\ cut along the
positive real axis.
\end{df}

This means that
$$
s_g =\mathcal{Q}(1_g)
$$
and we call $s_g$ the grading shift.

\subsubsection{Notation}
\label{shiftnotation}
In practice the choice of logarithms means that in a diagonal form
$\rho(g)=diag(exp(2\pi i \nu_1(g) ), \dots, exp(2\pi i \nu_n(g))$ and
$0\leq \nu_i \leq 1$.

For the element $j$,  $\nu_i(j)=q_i$, also due to choice of logarithm
$\nu_i(g)=1-\nu_i(g^{-1})$ furthermore
$d_g=\sum_{i:\nu_i(g)= 0} 1-2q_i$ and so

\begin{eqnarray}
s^+_g&=&2\sum_{i:\nu_i(g)\neq 0}(\frac{1}{2}-q_i),
\quad s^-_g= 2\sum_{i:\nu_i(g)\neq0} (\nu_i(g) -\frac{1}{2})\\
 s_g&=& \sum_{i:\nu_i(g)\neq 0}(\nu_i(g) - q_i).
\end{eqnarray}

\subsubsection{The bi--grading}
There is a natural bi--grading for the $G$--Frobenius algebras of
the type $G\M_f$ which is given by
$$
(\mathcal{Q},\bar {Q})(a_g):= (\Q(a_g)+s_g, \Q(a_g)+\bar s_g) \quad
\text{ for } a_g \in A_g
$$
where we used the notation $\bar s_g:= \frac{1}{2}(s_g^+ -s_g^-)$.

In the notation above
\begin{equation}
\bar s_g=\sum_{i:\nu_i(g)\neq 0} (1-\nu_i(g)-q_i)
=\sum_{i:\nu_i(g^{-1})\neq 0}(\nu_i(g^{-1})-q_i)=s_{g^{-1}} = d-d_g-s_g
\end{equation}

\begin{rk}
This grading is physically motivated, \cite{IV}, and basically
means that the natural bi--degree of the so--called ``twist
field'' is $(s_g,\bar s_g)$.
\end{rk}

\begin{lm}
An element $a_g \in A_g$ has diagonal grading $(q,q)$ if and
only if $s_g=s_{g^{-1}}$, i.e. $s_g^-=0$ or
equivalently $\sum_i \nu_i =\frac{1}{2} \codim(Fix(g))$.
\end{lm}

\subsubsection{Super--grading}
Recall that $A_g$ are all cyclic $A_e$ modules and the natural
parity for all of $A_e$ is all even.  Thus under the assumption
that all elements of $A_e$ are even the possible super--gradings
for the $\M_f$  module $G\M_f$ are given by maps $\sim \in
Map(G,\mathbb{Z}/2\mathbb{Z})$. Here $\tilde 1_g = \tilde g$.
Here we use $\sim$ both for the grading on $A$ and $G$ which
is justified by the equation above.

\subsubsection{The $G$--actions}
The primary choice for a $G$--action on $G\M_f$ would be the
induced action via pullback. Since $G$ acts on the collection of
fixed point sets: $h:\mathrm{Fix}_g \rightarrow
\mathrm{Fix}_{hgh^{-1}}$, we get a right action $r$ on $G\M_f$
which coincides with the notation of \cite{JKK}. On the other
hand, if we take the associated left action $l$, we are in line
with the definitions of \cite{orb} --- here $l(g):= r(g^{-1})$.

However, at this stage all actions of $G$ are good which preserve
the structures. Now  each $A_g$ is a cyclic $A_e$ module and we
denote the generator by $1_g$. We find that if $G$ is acting via
$A_e$ module automorphisms then

$$
\varphi(g) 1_h = \varphi_{g,h}1_{ghg^{-1}} \text { for some }
\varphi_{g,h} \in k^*
$$

From the fact that this is indeed an action of $G$, we obtain a
co--cycle condition on the $\varphi_{g,h}$. To be precise, they
form a non--abelian $G$ 2--co--cycle with values in $k^*$, where:

\begin{df} A non--abelian $G$ 2--cocycle with values in $k^{*}$ is a map
$\varphi: G\times G \rightarrow K^*$ which satisfies:
\begin{equation}
\label{grouphom} \varphi_{gh,k}=\varphi_{g,hkh^{-1}} \varphi_{h,k}
\end{equation}
 where $\varphi_{g,h}:= \varphi(g,h)$ and
 $$\varphi_{e,g}=\varphi_{g,e}=1$$
\end{df}

\subsubsection{The super--grading}
To define the super--grading we make an additional assumption.

{\bf Assumption} Keeping the condition that all elements of $A_e$
are even, we furthermore postulate that the $G$--action is an even
action.

This limits the possible super--gradings to functions of $C(G)
\rightarrow \mathbb{Z}/2\mathbb{Z}$ where $C(G)$ are the conjugacy
classes of $G$.

\subsubsection{The conditions from the trace axiom}
As demonstrated in \cite{orb}, if we further demand that the restricted trace
axiom holds\footnote{This means $[g,h]=e,c=1$} for the above $D(k[G])$--module, certain conditions for
the character, the super--grading and the co--cycle
$\varphi_{g,h}$ must hold.

Recall $\rho:G \rightarrow GL(n, {\bf C})$ is the representation
fixed from the beginning.

\subsubsection{The trace axiom, the character and the super--grading}

From the proof of the Theorem  5.1. \cite{orb} we extract the
following conditions proposition:

\begin{prop}
Let $G\M_f$ be the $D(k[G])$--module  with the $G$--action given by
the co--cycle $\varphi_{g,h}$ and fix a super--grading $\sim$ and a
character $\chi \in \mathrm{Hom}(G,k^*)$, then the trace axiom for
$g=e$ and arbitrary $h$ is satisfied w.r.t.\ $\sim$ and $\chi$ if
and only if $\chi$ satisfies
$$
\chi_{g}= (-1)^{\tilde g} (-1)^{n-\dim(\mathrm{Fix}_{g})}
\mathrm{det}(\rho(g))
$$
\end{prop}

\begin{rk}
Notice that this entails a condition on the super--grading:

Set
\begin{equation}
\label{sigmadef} \sigma(g):= \tilde g + n- \dim(\mathrm{Fix}(g))
\mod 2
\end{equation}
and call it the sign of $g$. Then
$$
\chi_{g}= (-1)^{\sigma(g)} \det(g)
$$
and therefore $\sigma \in \mathrm{Hom}(G,\mathbb{Z}/2\mathbb{Z})$.

Thus the possible super--gradings $\sim$ are in 1--1 correspondence
with elements $\sigma$ of $\mathrm{Hom}(G,\mathbb{Z}/2\mathbb{Z})$.
\end{rk}

\subsubsection{The trace axiom and discrete torsion}

\begin{df}
A discrete torsion bi--character for a group $G$ is a map from commuting pairs
$(g,h) \in G \times G: [g,h]=e$ to $K^*$ with the properties:
\begin{equation}
\label{epscond} \eps(g,h)=\eps(h^{-1},g) \quad \eps(g,g)=1 \quad
\eps(g_1g_2,h)= \eps(g_1,h)\eps(g_2,h)
\end{equation}
\end{df}

\begin{df} A non--abelian 2--cocycle  is said to satisfy the condition
of discrete torsion with respect to a given $\sigma \in {\rm
Hom}(G,\Zz)$ and a linear representation $\rho \in {\rm
Hom}(G,GL(n))$, if for all elements $g,h \in G: [g,h]=e$:
\begin{equation}
\label{eps}
\eps(g,h):=\varphi_{g,h}(-1)^{\sigma(g)\sigma(h)}
\det(g)\det(g^{-1}|_{{\rm Fix}(h)})
\end{equation}
is a discrete torsion.
\end{df}

\begin{rk}
 Due to the properties of $\varphi$ as a non--abelian
cocycle, $\varphi$ the second and third condition of discrete
torsion (\ref{epscond}) are automatically satisfied. If
furthermore $\g_{g,h}\neq 0$ then the first condition reduces to
$$
\det(g)det(g^{-1}|_{{\rm Fix}(h)})\det(h)\det(h^{-1}|_{{\rm
Fix}(g)})=1.
$$
\end{rk}

\subsubsection{The action of discrete torsion}
In \cite{disc}, we analyzed the phenomenon of
discrete torsion and
showed that the different choices $\varphi$
can by obtained from a fixed
$\dg$ module by tensoring with the twisted
group algebra $k^{\alpha}[G]$ with
$\a \in Z^2(G,k^*)$.

The corresponding discrete torsion
bi--character to such an $\a\in Z^2(G,k^*)$ is given
by $$\eps(g,h)=\frac{\a(g,h)}{\a(ghg^{-1},g)}$$
\cite{disc}

From the considerations of \cite{disc} one obtains:
\begin{lm}
For two choices of non--abelian cocycles $\varphi$ and $\varphi'$,
let $A(\varphi)$ and $A(\varphi')$ be the $\dg$
modules based on the k--module $G\M_f$, then
there is a group cocycle $\a\in Z^2(G,k^*)$, s.t.\
that $A(\varphi') \simeq A(\varphi)\otimes k^{\alpha}[G]$.
\end{lm}


From the Proof of Theorem 5.1. of \cite{orb} we can also extract the
following:

\begin{prop}
The $D(k[G])$--module $G\M_f$ with the $G$--action given by the
co--cycle $\varphi_{g,h}$ satisfies the super--trace axiom if and only if
there is a $\sigma \in{\rm Hom}(G,\Zz)$ s.t.\
 $\varphi_{g,h}$ satisfies the condition
of discrete torsion with respect to a  $\sigma \in {\rm
Hom}(G,\Zz)$  and the linear representation $\rho \in \mathrm
{Hom}(G,GL(n))$.
\end{prop}

\begin{crl}
If the group $G$ is Abelian then specifying a $G$ action
by a non--abelian cocycle $\varphi$ which satisfies the restricted trace axiom
for the resulting $\dg$ module is equivalent to specifying a
discrete torsion bi--character
$\eps(g,h)$ and a group homomorphism $\sigma \in \mathrm{Hom}(G,\Zz)$

\begin{equation}
\label{eps2}
\varphi_{g,h}= \eps(g,h)(-1)^{\sigma(g)\sigma(h)}\det(g^{-1})\det(g)|_{{\rm Fix}(h)}
\end{equation}
\end{crl}

\subsubsection{The $G$--Frobenius algebra structures} As explained in
\cite{orb} and \cite{sq} there is no fixed preferred $G$ Frobenius
structure on the $\M_f$ module above in general, but rather a set
depending on the choice of a so--called super--sign and a two
cocycle. The main result of \cite{orb} in this respect is:

\begin{thm}[\cite{orb}]
Given a natural $G$ action on a realization of a Jacobian
Frobenius algebra $(A_e,f)$ with a quasi--homogeneous function $f$
 let $A:=\bigoplus_{g\in G} A_g$ be as defined above up
to an isomorphism of Frobenius algebras on the $A_g$ then the
structures of super $G$--Frobenius algebra on $A$ are in 1--1
correspondence with triples $(\sigma,\g,\varphi)$ where $\sigma
\in {\rm Hom}(G, \Zz)$, $\g$ is a $G$--graded, section independent
cocycle compatible with the metric satisfying the condition of
supergrading with respect to the natural $G$ action, and $\varphi$ is
a non--abelian two cocycle with values in $K^{*}$ which satisfies
the condition of discrete torsion with respect to $\sigma$ and the
natural $G$ action, such that $(\g,\varphi)$ is a compatible pair.
\end{thm}

In many cases the equations for the co--cycles allow one to find a
unique multiplication up to the twist by discrete torsion.

We refer the reader to \cite{orb} for details.

The co--cycle $\g$ is a special type of $A_e$ valued group 2--cocycle
which defines on the cyclic generators multiplication via
$$
1_g \circ 1_h := \g(g,h) 1_{gh}
$$
the extra conditions ensure that the extension of this multiplication
using the cyclic $A_e$--module structures is well defined.
We usually write $\g_{g,h}$ for $\g(g,h)$

The function $\sigma$ is related to the super--sign as follows.
\begin{equation}
\sigma(g):= \tilde g +|N_g| \mod 2
\end{equation}
where  $|N_g|:= \mathrm{codim}(Fix(g))$ in $\mathbb{C}^n$.

Also:
\begin{df}
A cocycle $\g\in Z^2(G,A_e)$ is said to satisfy the condition of
supergrading with respect to a given a linear representation $\rho
\in {\rm Hom}(G,GL(n))$, if $\g_{g,h}=0$ unless $|N_h|+|N_g|+
|N_{gh}|\equiv 0 (2)$. Here $|N_g|:={\rm codim}({\mathrm
Fix}(\rho(g))$ is the codimension of the fixed point set of $g$.
\end{df}

\subsubsection{The metric}
The metric is constructed in two steps.

In the first step the metric is constructed from the metrics on
the individual Frobenius algebras $A_g:=\M_{f_g}$. First notice
that since $\mathrm{Fix}(g)=\mathrm{Fix}(g^{-1})$: $A_g=
A_{g^{-1}}$. Now $A_g$ has a non--degenerate pairing $\eta_g$
which we wish to view as a pairing $A_g \otimes
A_{g^{-1}}\rightarrow k$. We set
$$
\eta':= \bigoplus_{g\in G} \eta_g \in A^*\otimes A^*
$$
 In order to ensure the compatibility of the metric with the
multiplication and the twisted commutativity
$\g_{g,g^{-1}}=(-1)^{\tilde g}\varphi_{g,g^{-1}} \g_{g^{-1},g}$ we
need to rescale the metric:
$$\eta := \bigoplus_{g\in G} ((-1)^{\tilde g}\chi_g)^{1/2}\eta_g$$

For a discussion of the choice of the square root we refer to
\cite{orb}.

\subsection{The metric on the invariants}
As shown in \cite{orb} the $G$--invariants will be a Frobenius algebra
with respect to the metric $\eta$ if and only if $\chi(g)=\pm 1$.

\begin{rk}
In physics terms, this means that the spectral flow operator
$\mathcal U_{(1,1)}$ survives in the projection.
\end{rk}

\subsection{The dual of a quasi--Euler $G$--Frobenius algebra for a
quasi--homogenous singularity with  symmetries}

\subsubsection{The grading operator and the Euler condition}
Any non--trivial Frobenius algebra $\M_f$ stemming from a
quasi--homogenous singularity  has a non--trivial grading operator
$\Q$ as discussed above and
$$J:= exp(2\pi i \Q) = diag(exp(2\pi i q_1), \dots,\exp(2\pi i q_n))$$
generates a non--trivial finite cyclic group $\langle J \rangle
\subset GL(n)$ of order $ord(J)$ the order of $J$. Moreover fixing
$J$ as the generator we can identify this group with
$\mathbb{Z}/ord(J)\mathbb{Z}$ with a fixed generator $j$ acting
via $\rho(j)=J$.

\begin{rk}
This means that $G_{max}$ is non--trivial. Since any symmetry has
to preserve the grading $j$ it is in the center of $G_{max}$ and thus
any of the $G_{max}$ Frobenius algebras constructed in \S \ref{QMf} will be Euler. This
will also be true for any subgroup $H\subset G_{max}$ which
contains $\langle J \rangle.$
\end{rk}

\begin{lm}
\label{geuler}
If $\forall g: \eps(g,j^{-1})(-1)^{\sigma(g)(\sigma(j)+1)}=1$
then the corresponding $G$--Frobenius algebras
of \S \ref{QMf} will be $G$--Euler.
This is for instance the case if $\forall \eps(g,j) \equiv 1$
and $\sigma(j)=1$ or $\sigma \equiv 0$
\end{lm}

\begin{proof}
$$\varphi_{h^{-1}j,h}= \eps(hj,h)(-1)^{\sigma(h^{-1}j)\sigma(h)}
\exp(2 \pi i \sum_{i:\nu_i\neq0}(\nu_i - q_i))=
\exp(2 \pi i s_g)$$
\end{proof}

 {\bf Assumption:} We will assume that in the data $(f,G,\rho)$,
$\langle J \rangle \subset \rho(G)$, when considering duals
on the level of $\dg$--modules. Going to the algebra level
we postulate that the $G$--Frobenius algebra
structures above are Euler or quasi--Euler with fixed Eulerization.

\begin{rk}
The condition above holds for $G_{max}$, so if it does not hold
for a subgroup $H\subset G_{max}$, then if we are in a quasi--Euler
case with fixed Eulerization, we can
enlarge $H$ to $G_{max}$ and perform the dualization for the
$G_{max}$--Frobenius algebra and then reduce to the $H$--Frobenius
sub--algebra or the respective $D(k[H])$ module.
\end{rk}

\subsubsection{The dual $k$--module}

Given the $G$--Frobenius algebra $G\M_f$ its dual $k$--module is
defined as
$$
\check A_g = A_{gj^{-1}}= \M_{f|_{\mathrm{Fix}(gj^{-1})}}
$$
where $j$ is the group element defining the exponential of the
grading operator $\Q$ via $\rho(j)=\exp(2\pi i \Q)$.

\subsubsection{The dual $D(k[G])$--module} The $G$--module structure is
given by pulling back taction and scaling by $\chi$. In the
case of a singularity the character $\chi$ is determined by a
choice of sign function $\sigma \in \mathrm{Hom}(G,\Zz)$ given by
$\chi(g) =(-1)^{\sigma(g)}\det(g)$. If we denote
the $G$--action on $\check A$ by $\check \varphi$ then using the
$k$--module isomorphism $M: A_g \rightarrow A_{gj^{-1}}$
$$
\check {\varphi}(g) (\check a_h):= \chi(g) M \varphi(g) M^{-1} (\check a_h)
\in \check A_{hgh^{-1}}, \quad \text{for } \check a_h \in \check
A_g
$$
or if we denote $M(a)=: \check a$ and fix $\sigma \in
\mathrm{Hom}(G,\Zz)$ then for $\check a \in \check A_h$
$$
\check{\varphi}(g) (\check a):= (-1)^{\sigma (g)}\det(g)
\check{(\varphi(g)(a))} \in \check A_{ghg^{-1}}
$$

Using  equation (\ref{eps2}) for $\check a_h = M(a1_{hj^{-1}}) \subset \check A_h$

\begin{eqnarray}
\check \varphi(g)(\check a_h)&=&
\check \varphi(g) M(a 1_{hj^{-1}})\nn\\
&=&
\eps(g,hj^{-1})(-1)^{\sigma(g)(\sigma(h)+\sigma(j)+1)}
\det(g)|_{{\rm Fix}(hj^{-1})} M(a 1_{ghg^{-1}j^{-1}})
\end{eqnarray}

\begin{rk}
If  $\det(g)=(-1)^{\sigma(g)}$ then $\chi$ is trivial.
This means that the dual and the $G$--Frobenius algebra
have the same invariants.
\end{rk}


\begin{lm}
\label{ginv}
$\check 1_e$ is invariant if and only if $G M_f$ is
$G$--Euler.
\end{lm}

\begin{proof}
Since $Fix(j)=\emptyset$ as $f$ can be chosen to contain
no linear terms (the linear terms would actually only add
Eigenspaces of Eigenvalue one) the condition
$$
\check \varphi(g)(\check 1_e)=\check 1_e
$$ reads
\begin{equation}
\forall g \in G:\eps(g,j^{-1})(-1)^{\sigma(g)(\sigma(j)+1)}.
\end{equation}
This is precisely the condition to be $G$--Euler of Lemma \ref{geuler}.
\end{proof}

\begin{crl} Unless
$\forall g \in G:\eps(g,j^{-1})(-1)^{\sigma(g)(\sigma(j)+1)}
\det(g)|_{Fix(j)}=1$ there is no Frobenius structure on the
$G$ invariants of $\check GM_f$ for $GM_f$ with these invariants.
\end{crl}
\begin{proof}
Without this condition there will be no invariant unit for the
$(a,c)$ ring since for non--trivial grading $\check 1_e$ is
the only element of bi--degree $(0,0)$.
\end{proof}

{\bf Assumption:} Due to the content of the lemma, we will only consider
taking the invariants of a the dual of a fixed $\dg$ module structure
on $GM_f$ if it is
$G$--Euler.

\subsubsection{The bi--grading for the dual}
The bi--grading for the dual is given by the general formula (\ref{dualbigrading})
$$
\check {\mathcal{Q}}(\check a) = \Q(a)+\check s_g \quad \bar
{\check{\mathcal{Q}}} := \Q(a_g) +\bar{\check s}_g \quad \text{
for } \check a_g \in \check A_g
$$
In the notation of \ref{shiftnotation} this reads as\\
$\check s_g:=s_{gj^{-1}}-d=\sum_{i:\nu_i(gj^{-1})\neq 0}
(\nu_i(g)+(1-q_i)+\Theta(q_i-\nu_i)-q_i) -\sum_i 1-2q_i$\\
and $\check {\bar {s}}_g:=\bar s_{gj^{-1}}=
\sum_{i:\nu_i(gj^{-1})\neq 0}
(1-(\nu_i(g)+(1-q_i)+\Theta(q_i-\nu_i))-q_i)$ and thus

\begin{eqnarray}
\check s_g
&=& \sum_{\nu_i(gj^{-1})\neq 0} (\nu_i(g)+\Theta(q_i-\nu_i))- d_{gj^{-1}}\\
\bar {\check {s}}_g  &=&\sum_{\nu_i(gj^{-1})\neq
0}\Theta(q_i-\nu_i) -\nu_i(g)
\end{eqnarray}

\begin{rk}
An element in $\check A_g$ has anti--diagonal grading $(-q,q)$
if and only if $\Q(a_g)=\frac{1}{2}d_{gj^{-1}}$
\end{rk}

\subsubsection{The metric}
The metric is as in the general case the pulled back metric.
It will have group degree $j^2$
and will be homogenous of bi--degree $(-d,d)$.


\subsubsection{The degenerate $G$--Frobenius structure}
As remarked previously, for the dual $\dg$ module one
cannot expect a $G$--Frobenius structure, but what we
called a degenerate $G$--Frobenius algebra of group
degree $j$, which induces a $C(G)$
graded Frobenius structure on the invariants, in the sense of \cite{JKK}.

\subsection{Mirror symmetry for singularities}

In the framework of mirror symmetry a quasi--homogenous
function $f$ with an isolated singularity is considered
as a Landau--Ginzburg $B$--Model and hence has
as a $(c,c)$ ring the $(c,c)$ realization of $M_f$ and has a trivial
$(a,c)$ ring $A_1$.

\begin{df}
We call a $G$ Euler orbifold $G$--Frobenius algebra $A$,
together with a degenerate
$G$--Frobenius algebra of degree $j$ on $\check A$ a model for the
mirror dual
of a singularity $M_f$ if the $G$ invariants of $A$ are spanned by
$1\in A$ and the $G$--invariants of $\check A$ are the $(a,c)$ realization.

We also just say in short $(\check A)^G$ is the mirror dual to $M_f$.

We call two $G$ Euler orbifold $G$--Frobenius algebra $A$ and
an $H$ Euler orbifold $H$--Frobenius $B$,
together with a degenerate
$G$--Frobenius algebra of degree $j$ on $\check A$ and
 a degenerate
$H$--Frobenius algebra of degree $j'$ on $\check B$
 mirror dual pair if the $A^G=(\check B)^H$
 and $(\check A)^G=B^H$. In short we say, $A$ and $B$ are mirror dual.
\end{df}

Constructions for mirror pairs come form the \cred.

\section{A mirror theorem for simple singularities and other examples}
\label{examples}

In this section, we calculate the orbifolds and duals in several examples.
We will consider the first example in the greatest detail and then leave
slightly more details to the reader as we continue with
to make the text more concise.

The main result of this is the following
Theorem whose proof follows from the calculations below which
are collected in the table \ref{allcalcs}.

\begin{thm}
Let $f$ be one of the simple singularities $A_n,D_n,E_6,E_7$ and $E_8$
or a Pham singularity with coprime powers,
let $J$ be the exponential grading operator and $\Gamma:=\langle j \rangle$
with $\rho(j)=J$.
Then there is a projectively unique maximally non--degenerate
degenerate $\Gamma$--Frobenius algebra structure of degree $j$
on $\check {\Gamma M_f}$.
Moreover the invariants of the $G$--Frobenius algebra
$GM_f$ are one dimensional and yield the Frobenius algebra $A_1$,
while the invariants of the  $\check{\Gamma M_f}$ are
isomorphic as a bi--graded Frobenius algebra to $M_f^{(a,c)}$.

In short: the $A,D,E$ singularities  are
mirror self dual in the sense $(\check{\Gamma M_g})^{\Gamma}$
is the mirror dual.
\end{thm}

\begin{table}
\label{allcalcs}
$
\begin{array}{c|c|c|c|c|c}
M_f&\text{restriction}&G&\sigma&GM_f^G&(\check GM_f)^G\\
\hline
A_n&&\Znn&0&A_1&A_n\\
A_{2n-1}&&\Znn&1&A_1&B_n\\
A_{2n-1}&&\Zz&0&B_n&I_2(4)\\
A_{2n-1}&n\text{ odd for dual}&\Zz&1&D_{n+1}&A_1\\
A_{2n-1}&&\Zn&0&I_2(4)&B_n\\
D_{n+1}&&\Zzn&0&A_1&A_{2n-1}\\
D_{n+1}&n \text{ even}&\Zn&0&I_2(4)&B_n\\
D_{n+1}&n \text{ odd}&\Zn&0&A_1&D_{n+1}\\
D_{n+1}&&\Zz&0&B_n&I_2(4)\\
D_{n+1}&n\text{ odd for dual}&\Zz&1&A_{2n-1}&I_2(4)\\
A_{k_1-1}\otimes \dots \otimes A_{k_n-1}&
k_i \text { coprime}&
\Z/k_1\Z \times \dots \Z/k_n\Z&0
&A_1&A_{k_1-1}\otimes \dots \otimes A_{k_n-1}\\
E_6&&\Z/3\Z \times \Z/4\Z&0&A_1&E_6\\
E_7&&\znine&0&A_1&E_7\\
E_8&&\Z/3\Z \times \Z/5\Z&0&A_1&E_8\\
\end{array}
$
\caption{
Since all groups are cyclic $\eps\equiv 0$,
$\mathrm{Hom}(G,\Zz)=e$ or $\mathrm{Hom}(G,\Zz)=\Zz$
defining the entry in the column $\sigma$.
The conditions for the duals are the condition to be
be quasi Euler.
}
\end{table}

It would be tempting to conjecture that if
$\Gamma$ is the group  generated be the grading operator
 then $\check {\Gamma M_f}^{\Gamma}$
is the mirror dual. This is however not true as the example
of the elliptic singularity $P_8$ or a Pham singularity
of non--coprime exponents such as $x_0^5 + \dots +x_4^5 $ below show.
One obstruction
is that $\Gamma M_f^\Gamma$ is more than one--dimensional.
If $\Gamma M_f^\Gamma\simeq A_1$ we would however expect that
$\check {\Gamma M_f}^{\Gamma}$ is the dual.

\begin{conj}
For an insolated singularity $f$, let
$\Gamma$ be the group generated by the exponential grading operator $J$
if $\Gamma M_f^\Gamma \simeq A_1$ then $\check {\Gamma M_f}^{\Gamma}$
is the mirror dual.
\end{conj}

Also from the explicit calculations below, we obtain that.

\begin{thm}
The \credo holds and produces mirror pairs  for the cases self--dual cases with
the group $G=G_{max}=\Gamma$ the group generated
by the exponential grading operator and $H=e$ and also for
the listed
in table \ref{mirrorpairs}
\end{thm}

\begin{table}
\label{mirrorpairs}
$
\begin{array}{c|c|c|c|c|c}
T&G&H&K=G/H&((T/H)^H,(\check {T/H})^H)&
((T/H)/(K))^{K},(\check {(T/H)/K})^{K})\\
\hline
\begin{array}{c}
A_{2n-1}\\n \text{ odd}
\end{array}
&\Zzn&\Zz&\Zn
&(D_{n+1},A_1)&(A_1,D_{n+1})\\
A_{2n-1}&\Zzn&\Zz&\Zn&(B_n,I_2(4))&(I_2(4),B_n)\\
\begin{array}{c}
D_{n+1}\\n \text{ even}
\end{array}
&\Zzn&\Zz&\Zn&(B_n,I_2(4))
&(I_2(4),B_n)\\
E_6&\Z/3\Z\times \Z/4\Z&e \times \Zz&\Z/3\Z\times \Zz&(F_4,I_2(4))
&(I_2(4),F_4)\\
\end{array}
$
\caption{Mirror pairs from \cred\protect\footnotemark }
\end{table}
\footnotetext{For the first pair we choose $\sigma=1$ for $\Znn$
which restricts to $\sigma=1$ for $\Zz$ and $\sigma =0$ for $\Zn$
For the second and pair  we used the embedding of
$B_n$ into $A_{2n-1}$ as a subalgebra in the
untwisted sector and for the last
pair the embedding of $B_n$ into $D_{n+1}$  as a boundary
singularity to calculate the last column. For this we consider
the respective cyclic submodules over the sub--algebra corresponding
to $B_n$ in the untwisted sector.}

\begin{rk}
It is interesting to note that for $B_n$ in either
of its usual descriptions of folding $A_{2n-1}$ or $D_{n+1}$ we
obtain a non--trivial $(a,c)$ ring which is $I_2(4)$. The same holds true
for $F_4$ This
feature seems to distinguish  $B_n, F_4$ as simple {\em boundary} singularities.
\end{rk}

\begin{rk}

For a given Coxeter
group $W$ from the list
$A_n,B_n, D_n, E_6,E_7,E_8$, $H_3,H_4,F_4,G_2$ and $I_2(k)$
We denote by $W$ denotes the Frobenius algebra for the corresponding
Coxeter group. For the definition of the respective
Frobenius manifolds
we refer to \cite{dubrovinbook}.

Notice that we also use $B_2=I_2(4)$.
\end{rk}

\subsection{The case of $A_n$}

\label{ansection} The $A_n$ singularity is
given by the function $f:= x^{n+1}$. The maximal symmetry group is
given by $G:=G_{max}=\Znn$. Set $\zeta_n := \exp(2 \pi i
\frac{1}{n+1})$ then the exponential grading operator is
$J=\zeta_n id$ and $G= \langle j \rangle$ with $\rho(j)=J$.

\subsubsection{The $G$--graded $k$--module $G\M_f$}
Since
$$Fix_{j^i}= \left\{ \begin{tabular} {ll}$\mathbb{C}$& \mbox{if $i=0$}\\
$0$&\mbox{else}\\
\end{tabular}
\right.
$$
as $k$--modules $\M_f=A_n$ and $\M_{f|_0}=k = A_1$, where $A_k$ denotes the
Frobenius algebra of the $A_k$ singularity.

 The $k$--module $G\M_f$ is given by
$$
G\M_f:=\bigoplus_{i=0}^{n+1} A_{j^i}, \quad A_{j^i} :=
\begin{cases} A_n &i=0\\A_1 &i=1,\dots,n+1\end{cases},\quad G\M_f= A_n \oplus A_1
\oplus \dots \oplus A_1
$$
We denote the generator of the $j^i$ twisted sector $A_{j^i}$ by
$1_{j^i}$ and have the representation

$$
\rho(j^k)=\zeta_{n+1}^i = \exp(2\pi i \frac{k}{n+1})=\exp(2\pi i
\nu(j^k))
  \quad \nu(j^k)=\frac{k}{n+1}
$$

\begin{rk}
\label{annotation}
We wish to point out that the notation $A_g$ for $g\in G$ is the
standard notation for the twisted sectors of a $G$-Frobenius
algebra $A$, while the notation $A_n$ with $n\in \mathbb{N}$ is
the standard notation for the Frobenius algebra of the $A_n$
singularity. We hope that this will not lead to confusion, since one index is
a natural number and the other index an element in a finite group.
\end{rk}

\subsubsection{The grading}

The grading is determined by the shifts
$$
s^+_{j^i}=\begin{cases}0&\text{ if } i=0\\
\frac{n-1}{n+1}& \text{ if } i\in \{1,\dots,n\}\end{cases} , \quad
s^-_{j^i}=\begin{cases}0&\text{ if } i=0\\
 \frac{2i-n-1}{n+1}& \text{ if } i\in \{1,\dots,n\}\end{cases}
$$
Thus:
$$s_{j^i}=
\begin{cases}
0&\text { if } i=0\\
\frac{i-1}{n+1} & \text { if } i\neq 0
\end{cases}
$$

\subsubsection{The super-grading}
The choices of super--gradings $\tilde{}$
are determined by a choice of $\sigma \in \mathrm{Hom}(\Znn,\Zz)$ by
$$
\tilde 1_g \equiv |N_g| +\sigma(g)
$$
where

$$|N_{j^i}|=\begin{cases} 0 & \text { if } i=0\\
1 & \text { if } 1\leq i \leq n
\end{cases}
$$

Now if $n=2m$ is even, $\mathrm{Hom}(\mathbb{Z}/(2m+1)\mathbb{Z},\Zz)=e$,
so $\sigma(j^i)\equiv 0$.

If $n=2m-1$ is odd, $\mathrm{Hom}(\mathbb{Z}/(2m)\mathbb{Z},\Zz)=\Zz$
and there are two choices
for $\sigma$, either $\sigma(j^i) \equiv 0$ or $\sigma(j^i)F\equiv i \mod 2$.

Let us fix $\sigma \in \mathrm{Hom}(\Znn,\Zz)$.

\subsubsection{The bi--grading}
$$
(\mathcal{Q}, \bar {\mathcal{Q}})(1_{j^i})= \left\{
\begin{tabular} {ll}
$(0,0)$& for $i=0$\\
$(\frac{i-1}{n+1}, \frac{n-i}{n+1})$&else
\end{tabular}
\right.
$$

\begin{rk}
The elements with a $(q,q)$ grading are the elements $1,z, \dots ,z^{n-1}\in A_e$
and $1_{\frac{n+1}{2}}$ in the case that $n=2m-1$ is odd. The
latter element is not invariant under the whole group $\Z/(2m)\Z$, but is is
invariant under the subgroup $\Zz$ as we will show in \S \ref{dnexample} below.
\end{rk}

\subsubsection{The $G$--action}
We have already fixed $\sigma \in \mathrm{Hom}(\Znn,\Zz)$.

Since $\Znn$ is abelian, its action is determined by a choice of
discrete torsion $\eps$ by the trace axiom
$$
\eps(g,h):=\varphi_{g,h}(-1)^{\sigma(g)\sigma(h)}
\det(g)\det(g^{-1}|_{{\rm Fix}(h)})
$$
now since $\eps(j^i,j^k)=\eps(j,j)^{i+k}=1^{i+k}=1$, we find that
$\eps\equiv 1$ and this implies that

$$\varphi_{j^i,j^k}= \begin{cases} 1 &\text { if } k=0\\
(-1)^{\sigma(j^i)\sigma(j^k)}\zeta^{-i} &\text{ if } k \in
\{1,\dots,n\}\end{cases}$$

\subsubsection{The metric}
 After the necessary re--scaling, the metric is given
by
$$\eta(z^i, z^k)= \delta_{i+k,n-1} \quad \eta(1_{j^i},1_{j^k})=
\delta_{i+k,n+1} ((-1)^{\tilde j \;i}\zeta)^{i/2}z^{n-1}$$

\subsubsection{The $G$--Frobenius structure}
Using the reconstruction Theorem we have to find a cocycle $\g$
compatible with the
action defined by $\varphi$ above and the grading.

From the general considerations we know $\g_{j^i,j^{n-1-i}}\in A_e$
and $\deg(\g_{j^i,j^{n-1-i}}) = d-d_{j^i}= \frac{n-1}{n+1}$ which yields
$$
\g_{j^i,j^{n-1-i}}=( (-1)^{\tilde j i}\zeta^{i})^{1/2}\rho = ( (-1)^{\tilde j}\zeta)^{i/2} z^{n-1}
$$
for the other $\gamma$ notice that $\deg(1_{j^i}) +\deg (1_{j^k})= \frac{i+k-2}{n+1}$ while
$\deg(1_{j^{i+j}})= \frac{i+k-1}{n+1}$ if $i+k\neq n+1$, but there is no element
of degree $\frac{1}{n+1}$ in $A_{j^{i+k}}$ for $ i+k\neq n+1$, so that the respective
multiplication must yield zero if the condition is not met.

Hence
$$
\gamma_{j^i,j^k}=
\left\{ \begin{tabular} {ll}
$((-1)^{\tilde j}\zeta)^{i/2}z^{n-1}$& for $i+k=n+1$\\
$0$&else
\end{tabular}
\right.
$$

\subsubsection{The $G$--invariants}
Regardless of the choice of $\sigma$ the only invariant of $\Znn \M_{z^{n+1}}$ is
the identity $1\in A_e$.

\begin{prop}
The $\Znn$ invariants of the $\dg$-module $\Znn \M_{z^{n+1}}$ are spanned by $1\in Ae$ and
thus any $\Znn$--Frobenius algebras built on $\Znn \M_{z^{n+1}}$ has as its invariants
the Frobenius algebra $A_1$.
\end{prop}

\begin{rk}
This is the expected result since the dual of $(A_n,A_1)$ is $(A_1,A_n)$
according to \cite{V,IV}.
\end{rk}

\subsubsection{The dual $G$--graded $k$--module}

The dual $G$--graded $k$--module $\check GM_f:= \bigoplus_{g \in
\Znn} \check A_g$ is given by

$$
\check A_{j^{-i}} :=
\begin{cases}A_1 &i\in \{0,\dots,n-1\} \\A_n &i=n\end{cases}
\quad \check {G\Q}_f= A_1\oplus \dots \oplus A_1 \oplus A_n
$$
here again the remark \ref{annotation} applies.

\begin{rk} Notice that it is convenient to choose the generator $j^{-1}$
for the group
$\Znn$ instead of $j$.
\end{rk}

\subsubsection{The dual $\dg$ module}

The $G$--action on the $\check A_{j^{-k}}$ is given by:
$$
\check\varphi_{j^{-i},j^{-k}}= \varphi_{j^{-i},j^{-(k+1)}}\chi(j^{-i})=
\left\{
\begin{tabular} {ll}
$(-1)^{\sigma(j^{-i})}\zeta^{-i}$& for $k=n\equiv -1 \mod (n+1)$\\
$(-1)^{\sigma(j^{-i}) (\sigma(j^{-(k+1)})+1)}$& else
\end{tabular}
\right.
$$

\subsubsection{The bi--grading}
The bi--grading is given by
$$
(\check{\mathcal{Q}}, \bar{\check {\mathcal{Q}}})(\check1_{j^{-k}})=
\left\{ \begin{tabular} {ll}
$(-\frac{k}{n+1}, \frac{k}{n+1})$&$k\in\{0,1,\dots , n-1\}$\\
$(-\frac{n-1}{n+1},0)$& for $k=n$\\
\end{tabular}
\right.
$$


\begin{prop}
In the case that $\sigma \equiv 0$ the $\Znn$--invariants of
$\check \Znn \M_{z^{n+1}}=(\Znn A_n)^{\vee}$ are the linear subspace
$$\langle \check 1_e, \dots, \check 1_{j^{-(n-1)}}\rangle$$
This subspace is isomorphic as a graded $k$--module to $A_n$.

In the case that $n=2m-1$ is odd and $\sigma(j^i)\equiv i \mod 2$
the $\Znn$ invariant subspace is
$$\langle \check 1_e, \check 1_{j^{-2}},\dots, \check 1_{j^{-2m}}\rangle$$
This subspace is isomorphic as a graded $k$--module to the
$(a,c)$ realization of the sub--$k$--module $B_{m}\subset A_{2m-1}$
\end{prop}

\subsubsection{The metric on the dual $\dg$ algebra}
The metric is, after re--scaling the generators by a non--zero
factor given by the formulas
\begin{eqnarray*}
 \check \eta(1_{j^{-i}},1_{j^{-k}})&=&
\delta_{i+j,n-1} \text { for } i,j\in \{0,\dots,n-1\} \\
\check\eta(z^i 1_{j^{-n}},z^k 1_{j^{-n}})&=& \delta_{i+j,n-1}\\
\check \eta(1_{j^{-i}},z^k 1_{j^{-n}})&=&\check \eta(z^k
1_{j^{-n}},1_{j^{-i}})=0\text { for } i\in \{0,\dots,n-1\}
\end{eqnarray*}

\subsubsection{The degenerate $G$--Frobenius algebra structure}

There is a multiplication compatible with the bi--grading.
It is unique up to scaling of the generators and is
given by
\begin{eqnarray}
\check 1_{j^{-i}} \check 1_{j^{-k}} &=&
\begin{cases}
1_{j^{-(i+k)}}
& \text { if } i+k\leq n-1\\
0&\text { if } i+k\geq n
\end{cases}\\
\label{anramondmult}
\check 1_{j^{-i}} \check 1_{j^{-n}}&=&0 \text{ for } i \in \{0,\dots,n-1\}\\
z^{n-1-i}\check 1_{j^{-n}}z^{n-1-j}\check
1_{j^{-n}}&=&\begin{cases}
z^{n-1-(i+k)}1_{j^{-n}}&\text { if } i+k\leq n-1\\
0&\text { if } i+k\geq n\\
\end{cases}
\end{eqnarray}

The following statement is straightforward.

\begin{lm}
This multiplication renders the metric, invariant i.e.\ it satisfies
$\check \eta (\check a, \check b\check c )= \check \eta(\check a,
\check b\check c )$. Furthermore $\check \eta$ is the projectively
unique non--degenerate pairing compatible with the bi--grading and
the above multiplication is the projectively unique maximally
non--degenerate multiplication rendering the metric invariant.
\end{lm}

\begin{rk}
\label{anramond} The multiplication above is {\em not} compatible
with the grading and group grading. But changing the equation
(\ref{anramondmult}) to
\begin{equation}
z^{n-1-i}\check 1_{j^{-n}}z^{n-1-j}\check 1_{j^{-n}}=0
\end{equation}
yields a multiplication that is (a) compatible with the
bi--grading (b) compatible with the group grading and (c)
compatible with the $G$--module structure and thus is compatible
with the $\dg$--module. This multiplication does not, however,
render the pairing $\check \eta$ invariant.
\end{rk}

\begin{lm}
The multiplication of Remark \ref{anramond} is the projectively
unique maximally non--degenerate multiplication  which is
compatible with the $\dg$--module structure, i.e. turn the $\dg$
module into a $\dg$ module and co--module algebra.
\end{lm}

\begin{rk}
\label{anrmmetric} The metric $\check\eta'$ given by
$$
\check \eta'(\check a,\check b):= \begin{cases} \check\eta(\check
a,\check b)& \text { for }\check a\in A_{j^{-i}},\check b\in A_{j^{-k}}
 \quad
i,k \in
\{0,\dots,n-1\}\\
0& \text{else}\end{cases}
$$
is invariant w.r.t.\ the multiplication of Remark \ref{anramond}.
\end{rk}

\begin{rk}
The multiplication of Remark \ref{anramond} together with the
metric of Remark \ref{anrmmetric} which contains degenerate elements
and has a non--trivial annihilator of the whole algebra is reminiscent of
the appearance of so--called Ramond sectors in the theory of
spin--curves \cite{JKV,PV,P}. For a discussion, see \S \ref{spinparagraph} below.
\end{rk}

Collecting the results from above yields:

\begin{prop}
In the case $\sigma \equiv 0$, the $\Znn$ invariants of $\check
{\Znn \M_{z^{n+1}}}=\check {\Znn A_n}$ together with the
multiplication of Remark \ref{anramond} and the metric of Remark
\ref{anrmmetric} are a  degenerate $G$--Frobenius algebra
of charge $j$ which is a projectively
unique maximally non--degenerate algebra. The Frobenius algebra
given by the invariants with the grading $\cbQ$
is isomorphic to $A_n$ as graded Frobenius algebras.

 As bi--graded Frobenius
algebras, the  invariants of the dual are the
$(a,c)$ realization or $A_n$:$(\check {\Znn A^{(c,c)}_n})^{\Znn}= A_n^{(a,c)}$.

In this sense, $(\check {\Znn A_n})^{\Znn}= A_n$ is the mirror
and $A_n$ is mirror self--dual.

In the case $n=2m-1$ and $\sigma(j^i)\equiv i \mod 2$
$\check {\Znn A_n}$ together with the multiplication
of Remark \ref{anramond} and the metric of Remark \ref{anrmmetric}
are a degenerate $G$ Frobenius algebra of degree $j$
which is the projectively unique maximally
non--degenerate algebra. This Frobenius algebra of the group invariants
is isomorphic to the
Frobenius sub--algebra $B_{m}\subset A_n $ as graded Frobenius
algebras. In terms of the bi--grading the invariants dual of $\Znn A_n^{(c,c)}$ with
non--trivial $\sigma$ is $B_{m}^{a,c}$.
\end{prop}

\subsection{The case of $A_{2n-1}$ with a $\Zz$ action}

In the case of  $A_{2n-1}$, we can restrict ourselves to the
action of the subgroup $\Zz\subset \Zzn$, generated by $j^{n}=:-1$
which acts on $z$ by $z \mapsto -z$.

We now consider the singularity $A_{2n-1}$ with the group of
symmetries $\Zz$.

\subsubsection{The $G$--Frobenius algebras}
The data of the bi--graded $\dg$ module can be read of by
restricting the data of \ref{ansection}.

There is a unique twisted sector for the element $j^{n}$ and the algebra
$$
A_e=A_{2n-1}, \; A_{-1}= A_1=k \quad \Zz \M_{z^{2n}}=
A_{2n-1} \oplus k
$$

The $G$ action is again determined by the fact that $\Zz$ is
cyclic forcing $\eps\equiv 1$ and a choice of $\sigma \in
\mathrm{Hom}(\Zz,\Zz)$


$$\varphi_{-1,-1}= (-1)^{\sigma(-1)+1}$$

The bi--grading is given by
$$
s^+_{-1}=\frac{n-1}{n},\quad s^-_{-1}=0, \quad
s_{-1}=\bar s_{-1}=\frac{n-1}{2n}
$$

The super-grading is given by $\tilde 1_{-1}\equiv \sigma(-1) \mod 2$

The metric is
$$
\eta(z^i,z^k)= \delta_{i+k,2n-2}, \quad \eta(1_{-1},1_{-1})=1, \quad
\eta(z^i,1_{-1})=\eta(1_{-1},z^i)=0
$$

Taking into account the results of \cite{orb,sq}, there is a unique
$\Zz$--Frobenius algebra structure with the multiplication
$$
z^i\circ z^k=z^{i+k} \quad \text{ for }i+k \leq 2n-2,
z^i\circ z^k=0 \quad \text{for } i+k > 2n-2
$$
$$
 z^i\circ 1_{-1}=\delta_{i,0} 1_{j^{m-1}} \quad
1 \circ 1_{-1}=1_{-1}, \quad  1_{-1}\circ 1_{-1}=z^{2n-2}.
$$

This multiplication and the metric are compatible with the bi--grading and
yield a $\Zz$ Frobenius algebra in both cases.

Notice that since $det(g)=\pm 1$ the metric will make the
invariants into a Frobenius algebra.

In the case of $\sigma(j^m) \equiv 1 \mod 2$ we obtain as invariants
$$
\langle 1, z^2, \dots z^{2(n-1)}, 1_{j^{n}} \rangle
$$

The bi--grading of the invariants is  diagonal
and given by
$$((\frac{1}{n},\frac{1}{n}),\dots,
(\frac{n-1}{n},\frac{n-1}{n}),
(\frac{n-1}{2n},\frac{n-1}{2n})).$$

In the case  $\sigma (j^{m}) \equiv 0$  the space of invariants
is
$$
\langle 1, z^2, \dots z^{2(n-1)} \rangle
$$
and the multiplication, the metric and the bi--grading
is the restriction of the ones above.

\begin{prop} In total we obtain,
\begin{enumerate}
\item The $\Zz$ invariants of ${\Zz}\M_{z^{2m-2}}$ with the choice
$\sigma(j^n)\equiv 1 \mod 2$ are isomorphic as a bi--graded Frobenius
algebra to the $(c,c)$
model of $\M_{x^m+xy^2}=D_m$.

\item The $\Zz$ invariants of ${\Zz}\M_{z^{2m-1}}$
with the choice $\sigma(j^{n})\equiv 0 \mod 2$ are isomorphic as a
bi--graded Frobenius algebra to $B_{m}^{(c,c)}$.
\end{enumerate}
\end{prop}


\begin{rk}
The result above in which the invariants
of the untwisted sector of $A_{2n-1}$ yield $B_n$ is an instance
of what is called folding cf.\ \cite{zuber} and see
\S \ref{folding} below.
\end{rk}

\subsubsection{The dual $\dg$--module}
For the dual, we obtain two sectors

$$\check A_e= M(A_{j^{-1}})\simeq A_1,\quad \check A_{j^{n}}
=M(A_{j^{n-1}})\simeq A_1$$

The dual bi--grading is given by
$$
\check s_e=\bar{\check {s}}_e=0, \quad
\check s_{-1}=-\frac{1}{2} \bar{\check {s}}_{-1}=\frac{1}{2}
$$

\begin{rk}
In the case, that $\sigma\equiv 0$
or the case that $n$ is odd and $\sigma(j)=\sigma(-1)=-1$
the action is the restriction of the Euler $G$--Frobenius algebra
of \S \ref{ansection} and is thus quasi--Euler.

In both these cases the dual action is defined and is given by
$$
\check\varphi_{-1,1}= (-1)^{\sigma(-1)(\sigma(j)+1)}=1, \quad
\bar\varphi_{-1,-1}=(-1)^{\sigma(-1)(\sigma(-1)+\sigma(j)+1)}
=(-1)^{\sigma(-1)}
$$
\end{rk}


Since is $(A_{2n-1},\Zz)$ is not Euler, but only quasi--Euler,
we cannot pull back the metric, but due to the  grading there is
projectively only one compatible  homogenous metric.

\begin{prop}
Projectively there
is a  Frobenius algebra structure on the duals
compatible with the group grading which is isomorphic as a bi--graded
Frobenius algebra to the $(a,c)$ realization of the algebra
$I_2(4)$. In the case that $\sigma \equiv 0$
the invariants are $I_2(4)$ and in the case that $n$ is odd
and $\sigma(-1)=-1$ the invariants are $A_1$.
\end{prop}

\subsection{The case of $A_{2n-1}$ with symmetry group $\Zn$}

In the case of $A_{2n-1}$, we can also consider the symmetry group
$\Zn \subset \Zzn$ which is generated by $j^2$.

Again the group is cyclic and $\eps \equiv 1$. In the case that $n$ is even
there is only one possible choice of $\sigma \equiv 0$.
In the case that $n$ is odd there are two possible choices $\sigma \equiv 0$ or
$\sigma(j^{2k})\equiv k$. The later choice is not quasi--Euler, however.

\subsubsection{The $G$--Frobenius algebras}
The invariants can be read off from \S \ref{ansection}.
For $\sigma \equiv 0$ or for $\sigma(j^{2k})\equiv k$,
there are no invariants in the
twisted sector and the invariants in the untwisted sector are

$$
\langle 1, z^{n} \rangle
$$
The bi--degrees are $(0,0),(\frac{1}{2},\frac{1}{2})$.

\begin{prop}
Projectively there
is only one Frobenius algebra structure on these invariants
compatible with the group grading which is isomorphic as a bi--graded
Frobenius algebra to the $(c,c)$ realization of the algebra for the Coxeter group
$I_2(4)$. This is also the restriction of the respective multiplication
on the unique $\Zzn$ Frobenius algebra $\Zzn M_{z^{2n}}$.
\end{prop}

\subsubsection{The dual}

We can only consider the quasi--Euler choice $\sigma \equiv 0$.

The linear spaces for the dual $\dg$ module
are all one dimensional $\check A_{\gen^{2k}}=k \check 1_{\gen^{2k}}$
and are all invariant.

The bi--degrees are
$$
(\cQ,\cbQ) (\check 1_{\gen^{2k}})= (-\frac{k}{n}, \frac{k}{n})
(\check 1_{\gen^{2k}}) \quad k\in\{1,\dots, n-1\}
$$

\begin{prop}
The dual $\check \Zn A_{2n-1}$ affords a projectively unique  graded
$\Zn$--Frobenius algebra structure with
trivial $\Zn$ action  which is equal to the
Frobenius algebra of its invariants and is isomorphic to
the $(a,c)$ realization of $B_n$
\end{prop}

\subsection{The case of $A_{p-1} \otimes A_{q-1}$ especially $E_6$ and $E_8$}
We will consider the tensor product
$A_{p-1} \otimes A_{q-1}$  for coprime  $p,q$.

The corresponding quasi--homogenous singularity is given
by $f=x^p+y^q$. For these singularities $q_x=\frac{1}{p},
q_y=\frac{1}{q}, d=\frac{2(pq-p-q)}{pq}, \mu = p+q$

Since $p$ and $q$ are coprime $G_{max}=\Z/(pq\Z)=\Z/pZ \times \Z/q\Z$ which
is generated by the grading operator $\gen = (\gen_p,\gen_q)$
in the tensor representation for the symmetry groups of the
$A_n$ factors.

$$
\rho(\gen)=\left(\begin{matrix}\zeta_p&0\\
0&\zeta_q\end{matrix}\right)
$$
where as usual $\zeta_p=exp(2 \pi i \frac{1}{p})$ and
$\zeta_q=exp(2 \pi i \frac{1}{q})$.

\subsubsection{The $G$--Frobenius algebras}
$$
f_{\gen^i} =\begin{cases} z^p+z^q &\text{ if } i=0\\
z^p &\text { if } i = rp, r\in\{1,\dots,q-1\}\\
z^q &\text { if } i = rq, r\in\{1,\dots,p-1\}\\
0 &else
\end{cases}
$$
in these cases, we get the twisted sectors linearly isomorphic
to $A_{p-1} \otimes A_{q-1}, A_{p-1}, A_{q-1}$ and $A_1$.
The group is cyclic and
hence $\eps \equiv 0$. In the case that $pq$ is odd,
we only have the trivial choice $\sigma \equiv 0$.
In the case that it is even we also have the possibility to
set $\sigma(\gen^i)\equiv i$. We will leave the latter case to the
reader.

The action is given by
$$
\varphi_{\gen^k,\gen^i} =\begin{cases} 1 &\text{ if } i=0\\
\zeta_q^{-k} &\text { if } i = rp, \; r\in\{1,\dots,q-1\}\\
\zeta_p^{-k} &\text { if } i = rq, \;  r\in\{1,\dots,p-1\}\\
\zeta_q^{-k} \zeta_p^{-k} &else
\end{cases}
$$
and we see that only $1 \in A_e$ is invariant.

The grading is given by
$$
s_{\gen^i}=
\begin{cases} 0 &\text{ if } i=0\\
\frac{j-1}{q} &\text { if } i = rp=kq+j, \; r\in\{1,\dots,q-1\}\\
\frac{j -1}{p} &\text { if } i = rq=kp+j, \;  r\in\{1,\dots,p-1\}\\
 \frac{j-1}{p}+\frac{l-1}{q}&\text {if }  i = rp+j=kq+l
\end{cases}
$$
$$
\bar s_{\gen^i}=
\begin{cases} 0 &\text{ if } i=0\\
1-\frac{j+1}{q} &\text { if } i = rp=kq+j, \; r\in\{1,\dots,q-1\}\\
1-\frac{j + 1}{p} &\text { if } i = rq=kp+j, \;  r\in\{1,\dots,p-1\}\\
2- \frac{j+1}{p}-\frac{l+1}{q}&\text {if }  i = rp+j=kq+l
\end{cases}
$$

\begin{prop}
The $\Z/pq Z$ invariants of the unique $\dg$ module $\Z/pq Z A_{p-1} \otimes
A_{q-1}$ are one dimensional and are thus isomorphic to the Frobenius algebra
$A_1$.
\end{prop}

\subsubsection{The dual}

The dual action is given by
$$
\check \varphi_{\gen^k,\gen^i} =
\begin{cases} \zeta_q^{-k} \zeta_p^{-k} &\text{ if } i-1=0\\
\zeta_p^{k} &\text { if } i-1 = rp, \; r\in\{1,\dots,q-1\}\\
\zeta_q^{k} &\text { if } i-1 = rq, \;  r\in\{1,\dots,p-1\}\\
 &else
\end{cases}
$$
From this we obtain a $pq-p-q+1=(p-1)(q-1)$ dimensional space of invariants
spanned by
$$
\langle \check 1_{\gen^i} \rangle \quad  i-1 \equiv \hskip -10pt / \; 0
\mod p \text{ and } i-1 \equiv \hskip -10pt / \; 0
\mod q
$$
The grading is given by
$$
\check {s}_{\gen^i}=
\begin{cases}
0 &\text{ if } i=0\\
-\frac{2(pq-q-p)}{pq} &\text{ if } i=1\\
\frac{j-1}{q}+\frac{2}{p}-2 &\text { if } i-1 = rp=kq+j, \; r\in\{1,\dots,q-1\}\\
\frac{j-1}{p}\frac{2}{q}-2 &\text { if } i-1 = rq=kp+j, \;  r\in\{1,\dots,p-1\}\\
 \frac{j-1}{p}+\frac{l}{q}-2&\text {if }  i-1 = rp+j=kq+l
\end{cases}
$$

$$
\bar {\check {s}}_{\gen^i}=
\begin{cases}
0 &\text{ if } i=0\\
0 &\text{ if } i=1\\
1-\frac{j-1}{q} &\text { if } i-1 = rp=kq+j, \; r\in\{1,\dots,q-1\}\\
1-\frac{j-1}{p} &\text { if } i-1 = rq=kp+j, \;  r\in\{1,\dots,p-1\}\\
2- \frac{j-1}{p}-\frac{l-1}{q}&\text {if }  i-1 = rp+j=kq+l
\end{cases}
$$
where we choose $i \in \{0, \dots, pq-1\}$.

By comparing degrees we arrive at:

\begin{lm}
Let $i \equiv j \mod p, j \in {2, \dots p} $ and $i \equiv k \mod q$ ,
$k \in \{2,\dots,q\}$
then the map
$$\check 1_{\gen^i} \mapsto x^{p-j} y^{q-l}$$
induces an isomorphism of graded vector spaces
between the $\Z/pq\Z$
invariants of $(\Z/pq\Z M_{x^p+y^q})^{\vee}=(\Z/pq\Z A_{p-1}
\otimes A_{q-1})^{\vee}$
graded by $\cbQ$ and the graded Milnor ring $A_{p-1}\otimes A_{q-1}.$

Moreover as bi--graded space $((\Z/pq\Z A_{p-1}\otimes A_{q-1}))^{\vee})^{\Z/pq\Z}$
is the $(a,c)$ realization or the $A$--model of $A_{p-1} \otimes A_{q-1}$.
\end{lm}

By comparing the degrees and group degrees, one obtains:
\begin{prop}
There is a projectively unique maximally non--degenerate
degenerate $G$--Frobenius structure on
$(\Z/pq\Z A_{p-1}\otimes A_{q-1})^{\vee}$ whose invariants are the mirror dual
to $A_{p-1} \otimes A_{q-1}$
\end{prop}

\begin{crl}
If we restrict ourselves to the case $p=3,q=4$, we obtain the mirror
to the $E_6$ singularity  and for $p=3,q=5$ the mirror for the
$E_8$ singularity.

In the first case the invariants are
$$
\check 1_{\gen^i},\;  i \in \{0,2,3,6,8,11\} \quad
\text {corresponding to }
1,xy^2,y,y^2,x,xy
$$
with bi--degrees
$$
(0,0),(-\frac{5}{6},\frac{5}{6}),(- \frac{1}{4},\frac{1}{4}),
(-\frac{1}{2},\frac{1}{2}),(-\frac{1}{3},\frac{1}{3}),
(-\frac{7}{12},\frac{7}{12})
$$

In the second case the invariants are
$$
\check 1_{\gen^i}, \;  i \in \{0,2,3,5,8,9,12,14\}
\text {corresponding to }
1,xy^3,y^2,x,xy^2,y,y^3,xy
$$
with bi--degrees
$$
(0,0),(-\frac{14}{15},\frac{14}{15}),(-\frac{2}{5},\frac{2}{5}),
(-\frac{1}{3},\frac{1}{3}),(-\frac{11}{15},\frac{11}{15}),
(-\frac{1}{5},\frac{1}{5}),
(-\frac{3}{5},\frac{3}{5}),(-\frac{8}{15},\frac{8}{15})
$$

\end{crl}

\subsubsection{The case of $E_6$ and the relation to $F_4$}
Using the above calculations we can obtain a mirror pair for
$F_4$ from $E_6$
via the tensor product. For this we use that
$E_6=(A_2 \otimes A_3)$ and $G_{max}=\Z/3Z\times \Z/4\Z$.

\begin{prop}
$(F_4, I_2(4))$ and $(I_2(4),F_4)$ are a mirror dual pair obtained
from the \credo for $E_6$ with $G_{max}=\Z/3\Z \times \Z/4\Z$,
$H=e \times \Zz$ and $G/H= \Z/3\Z\times \Zz$.
\end{prop}
\begin{proof}
Using the group $\Zz$ acting via
$e\times \Zz$:
$E_6/(\Zz)=(A_2 \otimes A_3)/(e\times \Zz)= (A_2/e \otimes A_3/(\Zz))$.
Thus by the previous calculations:
$$
(((\Zz) E_6)^{\Zz},(((\Zz) E_6)^{\vee})^{\Zz})
=(A_2 \otimes I_2(4), A_1 \otimes I_2(4))= (F_4, I_2(4))$$

For the dual pair with $G_{max}=\Z/3\Z \times \Z/4\Z$, $H=e \times
\Zz$ and $G/H= \Z/3\Z\times \Zz$:
\begin{multline*}((A_2 \otimes
A_3/(\Zz))/(\Z/(3\Z) \times \Zz))^H= A_2/(\Z/3\Z) \otimes
((A_3/(\Zz))/(\Zz))^{\Zz}\\ = A_1 \otimes I_2(4)\simeq I_2(4)
\end{multline*}

\begin{multline*}
 (((A_2 \otimes A_3/(\Zz))/(\Z/(3\Z) \times \Zz))^{\vee})^{H}\\=
((A_2/(\Z/3\Z))^{\vee})^{\Z/3\Z} \otimes
(((A_3/(\Zz))/(\Zz))^{\vee})^{\Zz} = A_2 \otimes I_2(4))\simeq F_4
\end{multline*}

\end{proof}

\subsubsection{Certain Pham singularities}

The same reasoning holds true for the Pham singulaties of
coprime powers
$$
f=x_1^{k_1} + \dots x_n^{k_n} \text {with $k_i$ pairwise coprime}
= A_{k_1-1}\otimes A_{k_n-1}
$$
Let $\Gamma$ be the group generated by the grading operator
then $\Gamma =\Z/k_1\Z \times \dots\Z/k_n\Z$

\begin{prop}
The $\Gamma$ invariants of the $\Gamma M_f$ with the choice of
trivial $\sigma$ are $A_1$ and there is a degenerate maximally non--degenerate
multiplication on the dual $\check \Gamma M_f$ which is projectively unique
and the $\Gamma$ invariants of the dual are the $(a,c)$ realization of $M_f$.

In other words $\check \Gamma M_f$ is mirror for $M_f$ .
\end{prop}

\subsection{The case of $D_n$}
\label{dnexample}
Recall that $D_{n+1}=\M_{x^n+y^2z}$  with $q_1=q_x=\frac{1}{n},
q_2=q_y=\frac{n-1}{2n}, d=\frac{n-1}{n}$.
For $n>4$ the maximal symmetry group is
$G_{max}=\Z/2n\Z=\langle \gen \rangle$.
The maximal symmetry group in the case $n=4$ is larger
$G_{max}=\Z/3\Z \times \mathbb{S}_3$, but
also contains the group $\Z/6\Z$ generated by $\gen$.
We will make further comments about
the case $D_4$ below
in \S \ref{D4}.

If we fix $\zeta_{2n}:= \exp(2\pi i\frac{1}{2n})$, then

$$
\gen = \left( \begin{matrix} \zeta_{2n}^{2}&0\\
0&\zeta_{2n}^{-1}
\end{matrix}
\right)
\quad
\gen^i = \left( \begin{matrix} \zeta_{2n}^{2i}&0\\
0&\zeta_{2n}^{-i}
\end{matrix}
\right)
$$

$$
\gen^{n+1}= \left( \begin{matrix}
\exp(2\pi i \frac{1}{n})&0\\
0&\exp(2\pi i \frac{n-1}{2n})
\end{matrix}\right)= \exp(2\pi i \Q)=\rho(j)=J
$$

This implies that $G_{max}=\langle J \rangle$ if and only if $n$ is odd.

Since $e$ fixes both $x$ and $y$, $\gen^l$ fixes neither $x$ nor $y$ for
 $l\neq 0,n$ and $\gen^n$ fixes $x$ but not $y$,
we see that the orbifold data is as follows

\begin{tabular}{l|c|c|c|c|c|c|c|c|c}
$g \in \Zzn$&$f_g$&$\M_{f_g}$&$d_g$&$\nu_1(g)$&$\nu_2(g)$&$\frac{1}{2}s^+_g$&$\frac{1}{2}s^-_g$
&$s_g$&$\bar s_g$\\
\hline
$g=e=\gen^0$&$x^n+xy^2$&$D_{n+1}$&$\frac{n-1}{n}$&$0$&$0$&$0$&$0$&$0$&$0$\\
$g=\gen^l, 0<l<n$&$0$&$A_1$&$0$&$\frac{l}{n}$&$\frac{2n-l}{n}$&$\frac{n-1}{2n}$&$\frac{l}{2n}$
&$\frac{l+n-1}{2n}$
&$\frac{n-1-l}{2n}$\\
$g=\gen^n$&$x^n$&$A_{n-1}$&$\frac{n-2}{n}$&$0$&$\frac{1}{2}$&$\frac{1}{n}$&$0$&$\frac{1}{2n}$&
$\frac{1}{2n}$\\
$g=\gen^l, n<l<2n-1$&$0$&$A_1$&$0$&$\frac{l-n}{n}$&$\frac{2n-l}{2n}$&$\frac{n-1}{2n}$&$\frac{l-2n}{2n}$
&$\frac{l-n-1}{2n}$
&$\frac{3n-1-l}{2n}$\\
\end{tabular}

Comparing the degrees, we arrive at:
\begin{prop}
The only elements of bi--degree $(q,q)$ of the bi--graded $\dg$ module $\Zzn D_{n+1}$ are
the elements of the untwisted sector.
\end{prop}

\subsection{The $G$ action}
Since $\Zzn$ is cyclic there is only one choice of discrete
torsion which fixes the choice of $\varphi$ to be
$\varphi_{\gen^k,\gen^l}= (-1)^{\s(k)\s(l)}\det^{-1}(\gen^k)
\det(\gen^k|_{Fix_{\gen^l}})$
which reads

$$
\varphi_{\gen^k,\gen^l}= \begin{cases}
1&l=0\\
 (-1)^{\s(k)\s(l)}\zeta_{2n}^{-k}& l\neq 0,n\\
(-1)^{\s(k)\s(n)}\zeta_{2n}^k&l=n
\end{cases}
$$

\subsubsection{The metric}
The pairing on the twisted sectors is given by
$$
\eta(1_{\gen^k},1_{\gen^l})=\begin{cases}
\delta_{k+l,2n} \exp(2 \pi i \frac{k}{4n})&\text{ for } k \leq n\\
\delta_{k+l,2n} \exp(2 \pi i \frac{2n-k}{4n})&\text{ for } k \leq n
 \end{cases},
\quad \eta(x^k 1_{\gen^n},x^k 1_{\gen^n})= \delta_{k+l,n-2}\\
$$
all other pairings are zero except for the pairing on the untwisted sector,
which remains the pairing of $D_n$.

\begin{lm}
The $\Zzn$ invariant subspace is $\langle 1 \rangle$, thus
$$(D_{n+1}/\Zzn)^{\Zzn}=A_1$$ as a graded Frobenius algebra.
\end{lm}

\subsubsection{The $G$--Frobenius structure}

A straightforward calculation shows that:
\begin{prop}
There is projectively only  one $G$--Frobenius algebra structure compatible with the
bi--grading, which is given by
$$\g_{\gen^k,\gen^l}= \begin{cases} \delta_{k+l,2n}\exp(2 \pi i \frac{k}{4n})\text{ for } k \leq n\\
\delta_{k+l,2n}\exp(2 \pi i \frac{2n-k}{4n})&\text{ for } k \leq n\\
\delta_{l,n}1&\text{ for } k = n
\end{cases}
$$
\end{prop}

\subsubsection{The dual $\dg$ module}
The bi--grading is given by

$$
\check s_{\gen^{-k}}=\begin{cases}
-\frac{k}{2n}&\text { if } k \in \{0,1,2, \dots, n-2, n,\dots, 2n-1\}\\
-\frac{n-1}{n}&\text { if } k = n-1\\
\frac{3-2n}{2n}&\text { if } k = 2n-1\end{cases}
$$
$$
\bar{\check {s}}_{\gen^{-k}}=\begin{cases}
\frac{k}{2n}&\text { if } k \in \{0,1,2, \dots, n-2, n,\dots, 2n-1\}\\
0&\text { if } k = n-1\\
\frac{1}{2n}&\text { if } k = 2n-1\end{cases}.
$$

\begin{lm}
In the case that $n$ is even,
the only elements of bi--grading $(-q,q)$
of the dual $\dg$ module $\check {\Zzn D_{n+1}}$
are
$$
\langle \check 1_e, \check 1_{\gen^{-1}}, \dots, \check 1_{\gen^{-(n-2)}},
 y1_{\gen^{-(n-1)}},1_{\gen^{-n}}, \dots , 1_{\gen^{-(2n-2)}},
x^{\frac{n-2}{2}}1_{\gen^{-(2n-1)}}\rangle$$
and in the case that $n$ is odd, the elements of degree $(-q,q)$
 of the dual $\dg$ module $\check {\Zzn D_{n+1}}$
are
$$
\langle \check 1_e, \check 1_{\gen^{-1}}, \dots, \check 1_{\gen^{-(n-2)}},
x^{\frac{n-1}{2}}1_{\gen^{-(n-1)}}, y1_{\gen^{-(n-1)}},1_{\gen^{-n}},
\dots , 1_{\gen^{-(2n-2)}} \rangle.
$$
\end{lm}

\subsubsection{The dual $G$--action}
The dual $G$--action is given by

$$
\check \varphi_{\gen^{-k},\gen^{-l}}=
\begin{cases}
(-1)^{\sigma(\gen^k)(\sigma(\gen^l)+
\sigma(\gen^{n+1})+1)}&\text{ for } l \notin \{n-1,2n-1\}\\
(-1)^{\sigma(\gen^{k})}\zeta_{2n}^{k}&\text{ for } l = n-1\\
(-1)^{\sigma(\gen^k)(\sigma(\gen)+
\sigma(\gen^{n+1})+1)}\zeta_{n}^{k}&\text{ for } l= 2n-1
\end{cases}
$$

A longer but straightforward calculation shows

\begin{prop} For the different choices of $\sigma$ we obtain:
\begin{enumerate}

\item
In the case $\sigma\equiv 0$, the $\Zzn$ invariants of $(\Zzn D_{n+1})^{\vee}$ are
$$
\langle \check 1_e, \check 1_{\gen^{-1}}, \dots, \check 1_{\gen^{-(n-2)}},
y1_{\gen^{-(n-1)}},1_{\gen^{-n}}, \dots , 1_{\gen^{-(2n-2)}} \rangle.
$$
Their bi--degrees are
$\check{\mathcal{Q}}(\check 1_{\gen^{-k}})=-\frac{k}{2n}\check 1_{\gen^{-k}}$,
$\bar{\check{\mathcal{Q}}}(\check 1_{\gen^{-k}})=\frac{k}{2n}\check 1_{\gen^{-k}}$
for $k\in \{0,1,\dots,2n-1\}\setminus \{n-1\}$
and
$\check{\mathcal{Q}}(y\check 1_{\gen^{-(n-1)}})
=-\frac{n-1}{2n}y\check 1_{\gen^{-(n-1)}}$,
$\bar{\check{\mathcal{Q}}}(y\check 1_{\gen^{-(n-1)}})
=\frac{n-1}{2n}y\check 1_{\gen^{-(n-1)}}$.

This is the spectrum of the $(a,c)$ realization of $A_{2n}$ and
there is a projectively unique maximally non--degenerate $G$--Frobenius
structure on the dual whose invariants are the $(a,c)$ realization of
$A_{2n-1}$.

\item
In the case that $\sigma(\gen^k)\equiv k \mod 2$ and $n$ is even, then
$\sigma(j)\equiv 1$ and
the invariants are
$$
\langle \check 1_e, \check 1_{\gen^{-2}}, \dots , 1_{\gen^{-(2n-2)}},
x^{\frac{n-2}{2}}\check 1_{\gen^{-(2n-1)}}\rangle.
$$
Their bi--degrees
are $\check{\mathcal{Q}}(\check 1_{\gen^{-2k}})=-\frac{k}{n}\check 1_{\gen^{-2k}}$,
$\bar{\check{\mathcal{Q}}}(\check 1_{\gen^{-2k}})=\frac{k}{n}\check 1_{\gen^{-2k}}$
for $k\in \{0,1,\dots, n-1\}$
and
$\check{\mathcal{Q}}(x^{\frac{n-2}{2}}\check 1_{\gen^{2n-2}})
=-\frac{n-1}{2n}x^{\frac{n-2}{2}}\check 1_{\gen^{2n-2}}$,
$\bar{\check{\mathcal{Q}}}(x^{\frac{n-2}{2}}\check 1_{\gen^{-(2n-1)}})
=\frac{n-1}{2n}\check x^{\frac{n-2}{2}} 1_{\gen^{-(2n-1)}}$.
This is the spectrum of $D_{n+1}$. Furthermore there is a unique
maximally non--degenerate $\Zzn$--Frobenius algebra structure of charge $j$
on $\check \Zzn D_{n+1}$ which has as invariants the $(a,c)$--realization of
$D_{n+1}$.

So for $n$ even, $D_{n+1}$ is self--dual with the choice
of non--trivial $\sigma$.

\item In the case that $\sigma(\gen^k)\equiv k \mod 2$ and $n$ is odd, then
$\sigma(j)\equiv 0$ and
the invariants are
$$
\langle \check 1_{\gen^{-1}}, \check 1_{\gen^{-3}}, \dots , 1_{\gen^{-(2n-3)}},
x^{\frac{n-1}{2}}\check1_{\gen^{-(n-1)}}\rangle.
$$
Their bi--degrees
are $\check{\mathcal{Q}}(\check 1_{\gen^{-(2k+1)}})=-\frac{2k+1}{2n}\check 1_{\gen^{-(2k+1)}}$,
$\bar{\check{\mathcal{Q}}}(\check 1_{\gen^{-2k}})=\frac{2k+1}{n}\check 1_{\gen^{-(2k+1)}}$
for $k\in \{0,1,\dots, n-1\}$
and
$\check{\mathcal{Q}}(x^{\frac{n-1}{2}}\check 1_{\gen^{n-2}})
=-\frac{n-1}{2n}x^{\frac{n-1}{2}}\check 1_{\gen^{n-2}}$,
$\bar{\check{\mathcal{Q}}}(x^{\frac{n-2}{2}}\check 1_{\gen^{2n-2}})
=\frac{n-1}{2n}$.

This case is non $G$--Euler and we see that there is no Frobenius
algebra structure on the invariants, since the unit is missing.
This means that the prospective unit $\check 1_e$ is not invariant
and there is not even a degenerate $G$-Frobenius algebra structure on
the dual.
\end{enumerate}
\end{prop}

\subsection{The case of $D_{n+1}$ with the symmetry group $\Zn$}

Let $\Zn \subset \Zzn$ be the subgroup of even
powers $\Zn=\langle \gen^{2k} \rangle$.

\begin{rk}
Notice that this subgroup is Euler if and only if $n$ is odd. Also in
this case $G_{max}\neq \langle j \rangle$ and $\Zn \simeq \langle j \rangle$.
\end{rk}

Since most calculations are obtained via restriction from
those of the previous section, we handle both the $G$--Frobenius algebras
and the duals at the same time.

\subsubsection{The bi--gradings}

The calculations
above for the bi--grading for the $G$--Frobenius algebra and its dual
just restrict to sectors corresponding to the subgroup $\Zn$.

\subsubsection{The $\dg$--modules}

Since $\Zn$ is a cyclic group the discrete-torsion bi--character is trivial:
$\eps \equiv 1$.
In the case that $n$ is odd, there is only one choice $\sigma \equiv 0$.
In the case that $n$ is even, there are two choices for $\sigma$:
$\sigma \equiv 0$ and $\sigma(\gen^{2k}) \equiv k \mod 2$ and
In the first case the resulting structure is quasi--Euler,
while in the second case it is not.

\begin{prop}
For $n$ even and any choice of $\sigma$ the invariants
of the resulting $\dg$ module
on $\Zn D_{n+1}$ are two dimensional and are generated by
$\langle 1_e,x^{\frac{n}{2}}\rangle$. There is a projectively unique
Frobenius structure for the
invariants which is the structure of the Frobenius algebra $I_2(4)$.

If $n$ is odd, the invariants
of $\Zn D_{n+1}$ are one dimensional and generated by $1_e$.
Hence they are isomorphic to $A_1$ as Frobenius algebras.

In the case that $n$ is even, for the choice of $\sigma \equiv 0$
the resulting $\dg$ module structure
on the dual $\check {\Zn D_{n+1}}$ has invariants
$$
\langle \check 1_e, \check 1_{\gen^{-2}}, \dots , \check 1_{\gen^{-(2n-2)}}\rangle
$$
Their bi--grading is consistent with the $(a,c)$ realization of $B_n$
and the respective Frobenius algebra structure is compatible with the
group grading.
and
there is a degenerate $G$--Frobenius structure of charge $j$ on $\Zzn$
restriction to $\Zz$ is the dual $\check \Zn D_{n+1}$ and has
as invariants precisely the $(a,c)$ realization of $B_n$.

In the case that $n$ is odd and, the invariants are:
$$
\langle \check 1_e, \check 1_{\gen^{-2}}, \dots, \check 1_{\gen^{-(n-3)}},
x^{\frac{n-1}{2}}1_{\gen^{-(n-1)}}, y1_{\gen^{-(n-1)}},1_{\gen^{-(n+1)}},
\dots , 1_{\gen^{-(2n-2)}} \rangle.
$$
with bi--degrees matching the $(a,c)$ realization of $D_{n+1}$
and there a degenerate $G$--Frobenius structure of charge $j$ whose
invariants are precisely the $(a,c)$ realization of $D_{n+1}$.

So in the case of $n$ odd $D_{n+1}$ is mirror self dual
with  respect to the orbifolding by the symmetry
group generated by the grading operator
\end{prop}

\subsection{$D_n$ with the symmetry group $\Zz$}

In this subsection, we restrict the action of $G_{max}=\Zzn$ to the subgroup
$\Zz \subset \Zzn$ generated by $\gen^n=:-1$.

\subsubsection{The algebras $\Zz D_{n+1}$}

There are two twisted sectors which  as $k$--modules are
$$
A_e= D_{n+1} , A_{-1}= A_{n-1}
$$

$$
\phi_{-1,-1}=(-1)^{\sigma(-1)+1}
$$

There are two choices for $\sigma$, $\sigma(-1)\equiv 0$ or $\sigma(1)\equiv1$.
The first choice always yields a quasi--Euler $\Zz$ Frobenius algebra,
while the latter choice is quasi--Euler only in the case of $n$ odd.

The bi--grading and $\Zz$ action
can be read off from the tables in the previous section.

After fixing $\sigma$ there is a unique $\Zz$--Frobenius algebra
structure \cite{orb,sq} which is given by
$$
1_{-1} \circ 1_{-1}= x
$$

\subsubsection{The duals}
For the dual both $\check A_e$ and $\check A_{-1}$ are
one--dimensional and have degrees $(0,0),(-1/2,1/2)$.
Since this is at most a quasi--Euler
we cannot pull back the metric, but there is projectively
only one metric compatible with the group grading.

The action is given by

$$\check \varphi_{-1,1}=1 \quad
 \check \varphi_{-1,1}=(-1)^{\sigma(-1)}$$

\begin{prop}
In the case that  $\sigma(-1)\equiv 0$, the invariants are given by
$$
\langle 1,x,\dots,x^{n-1} \rangle
$$
The bi--grading and metric and multiplication are commensurate with those of
the $(c,c)$ realization of $B_{n}$.
The dual algebra has a projectively unique
Frobenius algebra structure compatible
with the bi--grading that is isomorphic to $I_2(4)$ and the invariants are
$A_1$.

In the case that $n$ is odd and $\sigma(-1)=-1$, the invariants
are
$$
\langle 1,x,\dots,x^{n-1},1_{-1},x1_{-1},\dots,x^{n-1}1_{-1} \rangle
$$
The algebra of invariants is isomorphic to the $(c,c)$
realization of $A_{2n-1}$ as
a bi--graded Frobenius algebra.

The dual algebra affords the structure of the $(a,c)$ realization of
$I_2(4)$ with trivial $\Zz$ action.
\end{prop}

\subsubsection{The case of $D_4$ and the relation to $G_2$}
\label{D4}

In the case $n=4$ the maximal symmetry group is
$\langle \gen, \frac{1}{2}\left( \begin{matrix} -1&i\\
-3i&1
\end{matrix}
\right)\rangle \subset GL(2,\mathbb{C})
$

Let $$j= \left( \begin{matrix} \zeta_3&0\\
0&\zeta_3
\end{matrix} \right)a=\left( \begin{matrix} 1&0\\
0&-1
\end{matrix}\right), b= \frac{1}{2}\left( \begin{matrix} -1&i\\
-3i&1
\end{matrix}\right)$$

Then $a^2=b^2=id, aba=bab$ and $\langle a,b \rangle \simeq \mathbb{S}_3$
the symmetric group on three elements. Also $\gen=aj$ and $\langle
\gen\rangle =\Z/6\Z=\Z/3\Z\times \Z/2\Z$.
Finally $G_{max}= \Z/3\Z\times \mathbb{S}_3$.

We do not want to present the full calculation, which is quite involved,
but note that the $G$--Frobenius algebra for $D_4/\mathbb{S}_3$
is given as a $k$--module by
$A_e=D_4, A_a=A_b=A_{aba}=A_2, A_{ab}=A_{ba}=A_1$. There are
three conjugacy classes and the invariants are
$1,x^2,1_{ab}\pm 1_{ba}$ where the sign is $+$ if
one uses $\sigma\equiv 0$ or $-$ if $\sigma(g)\equiv length(g)$.

For the invariance of $x^2$ notice that in the Milnor ring without
using an isomorphism $y^2=-3x^2$ and thus
$(\frac{1}{2}(-1+iy))^2=\frac{1}{4}x^2-\frac{1}{4}y^2+\frac{1}{2}ixy\equiv
\frac{1}{4}x^2-\frac{-3}{4}y^2=x^2$.

In the case of the group $\Z/3\Z$ generated by $ab$ there the $k$--module is
given by
$A_e=D_4,A_{ab}=A_{ba}=A_1$ and the invariants are
$1,x^2,1_{ab},1_{ba}$.

\begin{lm}
The invariants in the untwisted sector of $D_4/(\Z/3\Z)$ and
$D_4/\mathbb{S}_3$ are isomorphic to $G_2$ as graded Frobenius algebras.
\end{lm}

\subsection{The case $E_7$}

Recall that for $E_7: x^3+xy^3$, we have the following degrees
$q_1=q_x=\frac{1}{3},q_2=q_y=\frac{2}{9}, d=\frac{8}{9}$.

Fix $\zeta_{9}:= \exp(2\pi i \frac{1}{9})$ then the $E_7$
singularity has the exponential grading operator $J=\exp(2\pi
i \Q)$
$$
J=\left(
\begin{matrix}\zeta_{9}^{3}&0\\0&\zeta_{9}^{2}\end{matrix} \right)
$$
This operator generates a subgroup $\langle J \rangle \subset
GL(n,\mathbb{C})$ which is isomorphic to $\znine$. We fix a
generator $j$ of $\znine$ and regard the representation $\rho:
\znine \rightarrow GL(n,\mathbb{C})$ given by $\rho(j)=J$.

This is also the maximal symmetry group $G_{max}=\langle \gen \rangle$
$$
\gen=\left(
\begin{matrix}\zeta_{9}^{3}&0\\0&\zeta_{9}^{-1}\end{matrix} \right)
$$
and $J= \gen^7$.

\subsubsection{The $\znine$--graded $k$--module $\znine \M_f$}

The representation is given by
$$
\rho(j^i)=\left(
\begin{matrix}\zeta_{9}^{3i}&0\\0&\zeta_{9}^{2i}\end{matrix} \right)
$$

\begin{tabular}{c|c|r|r|r|r|r|r|r|r}
$g \in\znine$&$f_g$&$\M_{f_g}$&$d_g$&$\nu_1(g)$&$\nu_2(g)$
&$s^+_g$&$s^-_g$
&$s_g$&$\bar s_g$\\
\hline
$e=j^0$&$x^3+xy^3$&$E_7$&$\frac{8}{9}$&$0$&$0$&$0$&$0$&$0$&$0$\\
$j^1$&$0$&$A_1$&$0$&$\frac{1}{3}$&$\frac{2}{9}$&$\frac{8}{9}$&$-\frac{8}{9}$
&$0$&$\frac{8}{9}$\\
$j^2$&$0$&$A_1$&$0$&$\frac{2}{3}$&$\frac{4}{9}$&$\frac{8}{9}$&$\frac{2}{9}$
&$\frac{5}{9}$&$\frac{1}{3}$\\
$j^3$&$x^3$&$A_2$&$\frac{1}{3}$&$0$&$\frac{2}{3}$&$\frac{5}{9}$&$\frac{1}{3}$
&$\frac{4}{9}$&$\frac{1}{9}$\\
$j^4$&$0$&$A_1$&$0$&$\frac{1}{3}$&$\frac{8}{9}$&$\frac{8}{9}$&$\frac{4}{9}$
&$\frac{2}{3}$&$\frac{2}{9}$\\
$j^5$&$0$&$A_1$&$0$&$\frac{2}{3}$&$\frac{1}{9}$&$\frac{8}{9}$&$-\frac{4}{9}$
&$\frac{2}{9}$&$\frac{2}{3}$\\
$j^6$&$x^3$&$A_2$&$\frac{1}{3}$&$0$&$\frac{1}{3}$&$\frac{5}{9}$&$-\frac{1}{3}$
&$\frac{1}{9}$&$\frac{4}{9}$\\
$j^7$&$0$&$A_1$&$0$&$\frac{1}{3}$&$\frac{5}{9}$&$\frac{8}{9}$&$-\frac{2}{9}$
&$\frac{1}{3}$&$\frac{5}{9}$\\
$j^8$&$0$&$A_1$&$0$&$\frac{2}{3}$&$\frac{7}{9}$&$\frac{8}{9}$&$\frac{8}{9}$
&$\frac{8}{9}$&$0$\\
\end{tabular}

\begin{lm}
The elements of bi--degree $(q,q)$ of $\znine E_7$ are exactly the elements in the
untwisted sector $A_e$.
\end{lm}

\subsubsection{The $G$--action} For $\znine$ $\eps\equiv 1$ and $\sigma \equiv 0$,
so the $G$--action is given by
$$
\varphi_{j^i,j^k}= \begin{cases}1&\text {if } k=0\\
\zeta_9^{-2i}&\text {if } k\in \{3,6\}\\
\zeta_9^{-5i}&\text {else}
\end{cases}
$$
and the character is
$$
\chi(j^i)=\zeta_9^{5i}
$$

\begin{lm}
The $\znine$ invariants of the only compatible $D(k[\znine]$ module structure is given by
the unit $1_e$.
\end{lm}

\subsubsection{The dual bi--grading}
The dual grading is given by
$$
\begin{array}{c|c|c|c|c|c|c|c|c|c}
&0&1&2&3&4&5&6&7&8\\
\hline
\check s_{j^i}&0&-\frac{8}{9}&-\frac{8}{9}&-\frac{1}{3}&-\frac{4}{9}&-\frac{2}{9}&-\frac{2}{3}&
-\frac{7}{9}&-\frac{5}{9}\\
\bar{\check {s}}_{j^i}&0&0&\frac{8}{9}&\frac{1}{3}&\frac{1}{9}&\frac{2}{9}&\frac{2}{3}&\frac{4}{9}
&\frac{5}{9}
\end{array}
$$

The elements of bi-degree $(-q,q)$ are
$$
\langle \check 1_e, y^2 \check 1_{j}, \check 1_{j^2}, \check 1_{j^3},
 \check 1_{j^5}, \check 1_{j^6}, \check 1_{j^8} \rangle
$$
\subsubsection{The dual $\znine$ action}

The dual $\znine$ action is given by

$$\check \varphi_{j^i,j^k} =\begin{cases}
\zeta_9^{5i}&\text {if } k=1\\
\zeta_9^{3i}&\text {if } k\in \{4,7\}\\
0&\text{else}
\end{cases}
$$

\begin{lm} The $\znine$ invariants of the dual $\check \znine E_7$
are given by
$$
\langle \check 1_e, y^2 \check 1_{j}, \check 1_{j^2}, \check 1_{j^3},
 \check 1_{j^5}, \check 1_{j^6}, \check 1_{j^8} \rangle
$$
they are all of diagonal bi--degree, and their degrees are
$$
(0,0),(-\frac{4}{9},\frac{4}{9}),(\frac{-8}{9},\frac{8}{9}),
(-\frac{1}{3},\frac{1}{3}),( -\frac{2}{9},\frac{2}{9}),
(-\frac{2}{3},\frac{2}{3}),(-\frac{5}{9},\frac{5}{9})
$$
The pairing and bi--grading and the group grading
are commensurate with that of the
anti--chiral realization of $E_7$ under the association
$\check 1_e \mapsto 1, \check 1_{j}\mapsto y^2,
\check 1_{j^2} \mapsto x^2y, \check 1_{j^3}\mapsto x,
 \check 1_{j^5}\mapsto y, \check 1_{j^6}\mapsto x^2, \check 1_{j^8}\mapsto xy
 $, so that $E_7$ is self dual.
\end{lm}

Again by inspecting the grading and group grading
\begin{prop}
There is a unique maximally degenerate $G$--Frobenius structure
of charge $j$ on $\check \znine E_7$ whose invariants form
the $(a,c)$ realization of $E_7$. Hence $(\check \znine E_7)^{\znine}$
is the mirror dual to $E_7$.
\end{prop}

\subsection{The case $P_8$}
We would briefly digress to singularities of higher modularity.
The first singularity of this type is $P_8=x^3+y^3+z^3-axyz$ with
$a^3+27\neq 0$. The Milnor ring of this singularity is given
generated by $\langle 1,x,y,z,xy,yz,xz,xyz \rangle$.
It is quasi homogenous of degrees $q_x=q_y=q_z=\frac{1}{3}$ and
$d=1$.

This singularity is not self--dual. Moreover
in the case that $a\neq0$ there is no
symmetry group which has only $A_1$ as invariants of the
$G$--Frobenius algebra, since
the term $xyz$ always has to remain
invariant, so it is impossible for this a $G$--Frobenius algebra
built from this dual singularity to be mirror--dual
for any orbifolding group to another singularity. Also the invariants cease to
have the diagonal $(q,q)$ or anti--diagonal $(-q,q)$ grading.

Let us calculate $P_8/\Gamma$ for the group $\gamma$ generated
by the grading operator $J=diag(\zeta_3,\zeta_3,\zeta_3)$.
There are two one--dimensional twisted sectors.

The shifts for the twisted sectors $i=1,2$
are $s_{J}=0,\bar s_{J}=1; s_{J^2}=1,\bar s_{J^2}=0$
Since $det(J^i)=1$ and necessarily $\sigma\equiv0,\eps\equiv 1$
 all elements in the twisted sector are invariant.
 In total the invariant elements are
 $$
1,xyz, 1_J,1_{J^2} \text {of degrees } (0,0),(1,1),(1,0) \text{and } (0,1)
 $$

 For the dual, the action does not change since $\sigma \equiv 0$ and
 hence $\chi \equiv 1$
 and we obtain the same
 invariants, only with a shifted group grading.

\begin{rk}
Notice that the spectrum is such that is looks like the Hodge diamond
of manifold.
 \end{rk}

\begin{prop}
The $G$--Euler $G$--Frobenius algebra $\Z/3\Z P_8$ is mirror self--dual:
$(\Z/3\Z P_8)^{\Z/3Z}\simeq ((\Z/3\Z P_8)^{\vee})^{\Z/3Z}$.
\end{prop}
\section{Remarks on the relation to spin curves, the geometry
of singularities and folding}

\subsection{Remarks on the relation to $r$--spin curves and $A$--models for
quasi--homogenous polynomials}
\label{spinparagraph}

The $r$--spin curve picture was conceived by Witten as an
$A$--model or $\sigma$--Model counterpart for the $A_{r-1}$
Landau--Ginzburg $B$--model \cite{wittenspin}. In his construction
and the mathematical constructions of \cite{JKV,PV,P} this was
achieved. It turns out however, that in the formulation there are
two types of behaviors at given marked points called Ramond or
Neveu--Schwarz. The appearance of the Ramond case introduces an
additional element in the state space, which is $n+1$ dimensional
in the $A_n$ case. If this element appears in a correlation
function the value of the correlation function becomes zero. So
the algebra is what we called a degenerate Frobenius algebra of
degree $j$ if one assigns the group degree $j^{-1}$ to $z$ and
identifies the Ramond element with $z^{-n}$.

This is the projectively unique maximally
degenerate $G$--Frobenius algebra one obtains from $(pt/(\Znn))^{\vee}$
\cite{orb}. If one considers $A_1$ as
the $(a,c)$ ring of $A_n$ then by self duality of $A_n$
one could expect that $((A_1=pt)/\Znn)^{\vee})^{\Znn}=A_n$ (cf.\ \cite{orb})
which is indeed the structure found above.
In this interpretation the bi--grading is however not straightforward, although
the grading could be recovered from the $q_i$ and $\nu_i$ by considering
the action of $\Znn$ on $\mathbb{C}$ by roots of unity.

It would be desirable to consider not only this altered version
of our duality applied to the $(a,c)$ ring and not the $(c,c)$ ring,
but to see it directly on the $(c,c)$ side.

For this, we would like to give another interpretation of our previous
remark on the $A_{r-1}$ model. The equation here which is mimicked
in the case of spin curves is $z^r=0$ only that in
the spin--curve picture $\mathcal{L}^{\otimes r}\simeq\omega(twisted)$
(Here $\omega(twisted)$ is a suitably twisted version of the canonical line
bundle on the curve).
 The fact that the Ramond
sector is zero in all the correlation functions can be taken to mean
that it is in fact zero. In other words it appears only as a
degenerate state and as we project to the invariants the equation
$z^r=0$ is implemented. In algebraic terms for $A_{r-1}$, first one considers
$R:=\mathbb{C}/(z^r)$ and then $R/(z^{r-1})$. The first quotient
is inherent in the spin picture in the periodicity with respect to
$\omega(twisted) \simeq \mathcal{L}^{\otimes r}$.

Now our degenerate Ramond sector is in fact $n$ dimensional for $A_n$
and not one dimensional. Here one should remark that for the
construction of of an operad  in the Ramond case,
one would actually have to fix a choice of isomorphism of the line bundle
with $\omega(twisted)$. The space of choices for this isomorphism is
a principal $\Znn$ space and thus
if one includes this data in the moduli problem the
state space for the Ramond sector becomes  $n+1$ dimensional
for $A_n$\cite{Jprivate}.
So indeed the Ramond sector
seems to be intrinsically higher dimensional. The fact that the dimension
is not $n$, but $n+1$ dimensional can be understood by the reasoning before.
In our description
the singularity in this sector is the singularity $A_n$. In the Milnor ring
interpretation this produces a Frobenius algebra which has $n$ states.
In the spin--representation as discussed above one would expect
$n+1$ states, one of which is degenerate.

The musings on this subject are at the moment only on the level of the
un--deformed
algebra,
but we hope to make them into more solid statements.

There is a straightforward way to build a spin curve like picture
for any quasi--homogenous polynomial $f$.  For this one considers
a line bundle $\mathcal{L}_i$ for each of the variables $z_i$ and
imposes the equations obtained by substituting the line bundles
$\mathcal{L}_i$ into the of monomials of the polynomials $f$ and
equates these expression to $\omega(twisted)$. This defines the
moduli problem.
 This approach is being seriously discussed by
\cite{FJR}. When the polynomial $f$ is such that is maximal
symmetry group $G_{max}$ is Abelian and each variable appears by
itself, the corresponding virtual fundamental class is constructed
in \cite{FJR}. The hope is to be able to lift these conditions
\cite{Jprivate}. We would like to point out that the condition on
the variables appearing alone in a monomial ensures that
$\mathbb{C}[[z_i]]/(m_j)$ is finite dimensional. Here the $m_j$
are the monomials of $f$.

Further evidence for our interpretation of ``Landau--Ginzburg $A$--models''
arises from these constructions.
For each element $g\in G_{max}$ there are again two types of behaviors at
the marked point
which are either of Ramond or of Neveu--Schwarz type. The Ramond means that
the isotropy at a marked point is not the full symmetry group while in the
Neveu--Schwarz case it is.

Again to turn the resulting moduli spaces into operads it is necessary to include
additional data for the Ramond case which is isomorphic
to the reduced symmetry group of $f_g$ \cite{Jprivate}.

\begin{conj}
We conjecture that the Neveu--Schwarz sectors are in 1--1
correspondence with the one--dimensional twisted sectors and the
Ramond sectors are in one--one correspondence with the sectors
that are more than one--dimensional, i.e. $\Q_{f_g} \neq
\mathbb{C}$.
\end{conj}

This conjecture has been checked against the preliminary results
of \cite{FJR}.

\begin{conj}
We expect that the non--degenerate part of the cohomological Field
theory described by a quasi--homogenous polynomial is the
deformation of the Frobenius algebra of the invariants of
$(G_{max}M_f)^{\vee}$. Moreover, we expect that the behavior of
the correlation functions is modelled by the deformations of a
degenerate $G$--Frobenius of charge $j$ given by
$(G_{max}M_f)^{\vee}$, possibly adding more degenerate elements.
More precisely, let $m_j$ be the monomials of $f$ and $q_i =
\frac{1}{n_i}$ be the quasi--homogenous degrees of the $z_i$. In
the case that the ring $\hat M_{f_g} :=\mathcal{O}/(m_{f_g,j})$ is
finite dimensional, the extra elements should correspond to the
extension of basis from $M_{f_g}$ to $\hat M_{f_g}$ --- for each
higher--dimensional sector.
\end{conj}

Our calculations predict that this procedure yields the right result in
the case of Pham singularities with coprime powers, such as $E_6$ and
$E_8$ and indeed this is true by taking tensor products of spin--curves
\cite{JKV2}.

\subsection{Orbifolding and the geometry of singularities with symmetries}
\label{singtheory}

There is a relationship of our constructions of $G$--Frobenius
manifolds for a singularity $f$ and the Ramond state space of
\cite{orb} (not to be confused with the Ramond notation for
spin--curves) and classic singularity theory.

For a singularity there are classically two objects which are
studied, one is the Milnor ring $M_f$ which also provides
a basis for the minversal unfolding which can be written as
$$
F: (\mathbb{C}^{n+1}\times M_f,0) \rightarrow (\mathbb{C},0)
$$

This fact affords an extension by the choice of a primitive form \cite{Saito}
to a construction of Frobenius manifold on the flat space $M_f$
\cite{dubrovinbook}.

The other object of interest obtained from the Milnor fibration which
is given by
$$f:(\mathbb{C}^{n+1},0) \rightarrow (\mathbb{C},0)$$
which gives a local fibration on $\mathbb{C}-0$. The fibers are
bouquets of spheres and the Betti number of the middle dimensional
cohomology of the fibers is also $\mu$. The isomorphism between
$M_f$ and $H:= H^n(F_*,\mathbb{C})$ can be given by a choice of primitive
form. Here $F_*$ denotes a generic fiber.

Now suppose $G\subset GL(n+1,\mathbb{C})$ is a group of symmetries.
This will act on the total space of the Milnor fibration
and trivially on base and thus there is an induced action on $H$.

Let $det$ be the one dimensional representation of $G$ given by the determinant.

The main result of  \cite{Wa} is

\begin{thm}[\cite{Wa}]
In the situation described above the $\mathbb{C}[G]$ modules
$H$ and $M_f \otimes det$ are isomorphic.
\end{thm}

This infers that while the untwisted sector of the $G$--Frobenius algebra is
isomorphic as a $\mathbb{C}[G]$--module to $M_f$
 the untwisted sector of the Ramond state space
is isomorphic as a $\mathbb{C}[G]$--module to $H$.
This untwisted Ramond sector corresponds to the $j$ twisted sector of
the dual.

In  exactly the case that
 the symmetries generate a Coxeter group $G$ the quotient of
 $\mathbb{C}^{n+1}$ by $G$ is smooth:
 $\mathbb{C}^{n+1}/G \simeq \mathbb{C}^{n+1}$.
In this situation, one
can regard the germ $f_G$ on the quotient. Let $\mu_G$ denote
the Milnor number of $f_G$ and $\mu_g$ those of $f_g:=f|_{Fix(G)}$.
Here we need to assume that this restriction is again an isolated singularity
which is automatic in the quasi--homogenous case.
Also fix $d_g = codim(Fix(g))$

The results of \cite{Wa} are

$$
\mu_G=\frac{1}{|G|}\sum_{g\in G}(-1)^{d_g} \mu_g
$$

Furthermore in \cite{Wa} the equivariant Euler--characteristic of the
 $\mathbb{C}[G]$--modules $M_f$ and $H$ is used to compute the Milnor numbers
$\mu_g$. Let $M=H_n(F_*,\mathbb{C})$ and consider its class
$[M]$ in the representation ring of $G$. We can identify this with the ring
of class functions and evaluate at elements $g$.

The formula is \cite{Wa}
$$
\mu_g=(-1)^{d_g}[M](g)
$$

This give a way to compute the invariants of the untwisted Ramond state
space which is isomorphic as a $G$--module to the $j$ twisted sector
of $\check GM_f$.

It is interesting to note that the twisted sectors contribute to this
calculation through the equivariant Euler characteristic.

One could adapt these techniques to the restrictions of the singularity to
the various fixed point sets and obtain formulas for the dimension of
the whole space $\check M_f$.

\subsection{Folding}
\label{folding}
For the Dynkin diagrams of the simple singularities and
more generally for the generalized Dynkin diagrams of
\cite{zuber}, there is an operation known as folding.

In this section we show, that the folding can be described
as a non--stringy orbifolding with respect to a group of projective
symmetries.

\begin{df}
A projective symmetry for a singularity $f:\mathbb{C} \rightarrow \mathbb {C}$
with an isolated critical at zero is an element
$S\in
GL(n,\mathbb{C})$, s.t.\ $f(S(\mathbf{z}))=\lambda f(\mathbf{z})$ for some
$\lambda \in \mathbb{C}$

A projective folding group for a quasi--homogenous singularity $f$ is group
$G$ together with a representation of $G$ in $GL(n,\mathbb{C})$ which
acts by projective symmetries with the same fixed $\lambda$ and
preserves the unique (up to scalar multiples) element of highest degree.
\end{df}

These type of symmetries act on the Milnor ring, since
the local ring $f({\bf z})=0$ is equal to that of $\lambda f({\bf z})=0$.

\begin{rk}
For a sum of two singularities $f+g$, the product
of two projective symmetry groups for $f$ and $g$ respectively
also acts on the Milnor ring $M_{f+g}=M_f \otimes M_g$
\end{rk}

\begin{thm}
For each of the classical foldings for Coxeter groups there is
a group of projective symmetries or a product
of two groups of projective symmetries which has as its invariants the
Frobenius algebra of the folded graph. The foldings and groups
are contained in table \ref{foldingtable}.
\end{thm}

\begin{table}
\label{foldingtable}
$
\begin{array}{c|c|c|c}
\text{Diagram/group}&\text{Folded diagram/group}&
\text{Folding group}&\text{representation}\\
A_n&I_2(n+1)&\Z/(n-1)\Z&z \mapsto \zeta_{n-1}z\\
A_{2n-1}&B_n&\Zz&z\mapsto -z\\
D_{n+1}&B_n&\Zz&(x,y)\mapsto (x,-y)\\
D_4&G_2&\Zz^*&(x,y)\mapsto (-x,-y)\\
D_6&H_3&\Zz&(x,y)\mapsto (-x,-y)\\
E_6&F_4&\Zz&(x,y)\mapsto (x,-y)\\
E_8&H_4&e\times \Z/3\Z&(x,y)\mapsto (x,\zeta_3 y)
\end{array}
$


{\footnotesize ${}^*$ This is the simplest group. Other folding groups are
$\Z/3\Z$ and $\mathbb{S}_3$ as discussed in \S \ref{D4}.}
\caption{The Foldings and their projective symmetry groups.}

\end{table}

\begin{rk}
The utilization of projective symmetries is necessary, since not
all foldings can be realized
with $\lambda=1$ in particular the element of highest degree
transforms in the representation $\det(\rho(g))^{-2}$ (see e.g. \cite{orb})
so that the only folding groups with $\lambda =1 $ will be those whose
determinant lie in $\pm 1$. In particular the $\Zz$ foldings of
$A_{2n-1}$ and $D_{n+1}$ yielding as discussed above and also
$E_6$ to $F_4$ can be realized by orbifolding.
For $G_2$ the folding can be only be obtained via
orbifolding by restricting to the classical level,
i.e. disregarding the twisted sectors.
\end{rk}

\begin{rk}
The folding of $E_6$ and $E_8$ can also be understood as
the tensor products of the folding on the factors.
$A_2\otimes I_2(4)=F_4$ and $A_2\otimes I_2(5)$.
\end{rk}

\begin{rk}
Unlike in the case of the operation of symmetries, the group action
 of projective symmetries does
 not act on the Milnor fibration fiberwise and hence not obviously
 on the cohomology bundle. But on the other hand it leaves the central
 fiber invariant and furthermore acts by homothety on the base via
 $f(z)=t \mapsto f(z)=\frac{1}{\lambda}t$, so we obtain an equivariant
 action.
\end{rk}

\begin{rk}
The relation of folding to the miniversal unfolding space is known
 and is given in \cite{strachan}. In fact the foldings provide
 submanifolds of Frobenius manifolds or $\mathcal{F}$ manifolds.
\end{rk}

\begin{rk}
It would be desirable  to extend the theory of $G$--Frobenius algebras
to these quotients. In fact it seems to be straightforward to generalize
some of the construction of \cite{orb} for Jacobian Frobenius
algebras with symmetries
to those with projective symmetries. Here the twisted sectors would again
just be obtained from the function by restriction to the fixed subspace.
There is also no obstruction to keeping the grading shifts and dualization
process. One would expect to be able to apply this type of orbifolding to
the calculations and definition of \cite{zuber}.
 We leave the more careful analysis of this possibility for the
future.
\end{rk}


\begin{thebibliography} {99}

\bibitem
[AGLV]{arnoldbook}
 V.\ I.\ Arnold et al.
{\it Singularity theory. I.} Translated from the 1988 Russian original
by A. Iacob.
Reprint of the original English edition from the series Encyclopaedia
of Mathematical Sciences [ Dynamical systems. VI,
 Encyclopaedia Math. Sci., 6, Springer, Berlin, 1993

\bibitem[CR]{CR} W.~Chen and Y.~Ruan, {\it A new cohomology theory for
orbifold.} {math.AG/0004129}. And
W.~Chen and Y.~Ruan, {\it Orbifold Gromov-Witten theory.}
In A. Adem, J. Morava, and Y. Ruan
(eds.), {Orbifolds in Mathematics and Physics}, {\it Contemp.
Math.}, Amer. Math. Soc., Providence, RI. {310}, (2002),
25--85. {math.AG/0103156}.

\bibitem
 [Du]  {dubrovinbook}
B.\ Dubrovin,  Geometry of 2D topological field
theories. in {\it Integrable systems and quantum groups
(Montecatini Terme, 1993)}, Lecture Notes in Math., 1620,
Springer, Berlin, 1996, pp.\ 120--348.



\bibitem
[FJR]  {FJR} H.\ Fan, T.\ Jarvis and Y.\ Ruan, {\it A
Generalization of Spin Orbifold Quantum Cohomology Arising From
Quasi-Homogeneous Polynomials}. In preparation.


\bibitem
[GP]   {GP}
B.\ R.\ Greene,  M.\ R.\ Plesser,
Duality in Calabi-Yau moduli space,
{\it  Nuclear Phys. B} {\bf 338} (1990), 15--37.

\bibitem
[IV]{IV}
K.\ Intriligator and C.\ Vafa,
 Landau-Ginzburg orbifolds,
{\it Nuclear Phys. B} {\bf 339} (1990), 95--120

\bibitem
[J]{Jprivate} T.\ Jarvis. Private communication.

\bibitem
[JKK] {JKK}
T.\  Jarvis, R.\  Kaufmann and T.\ Kimura,
Pointed Admissible $G$-Covers and
$G$-equivariant Cohomological Field Theories, Preprint
math.AG/0302316.
14 (2003), 573-619.

\bibitem
[JKV]   {JKV}
T.\ Jarvis, T.\ Kimura and A.\ Vaintrob.
{\it Moduli spaces of higher spin curves and integrable hierarchies}
{Compositio Math.} {\bf 126} (2001).

\bibitem
[JKV2]   {JKV2}
T.\ Jarvis, T.\ Kimura and A.\ Vaintrob.
{\it Tensor products of Frobenius manifolds and moduli spaces of
higher spin curves.} Conf\'erence Mosh\'e Flato 1999, Vol. II (Dijon),
145--166, Math. Phys. Stud., 22, Kluwer Acad. Publ., Dordrecht, 2000.
\bibitem
[Ka1]   {tenfrob}
R.\ M.\ Kaufmann,  {\it
The tensor Product in the Theory of Frobenius manifolds},
{ Int.\ J.\ of Math.} { 10} (1999) 159--206.
\bibitem
 [Ka2] {conf}
R.~M.~Kaufmann,
{\it Orbifold Frobenius algebras, cobordisms,
and monodromies}. In A. Adem, J. Morava, and Y. Ruan (eds.),
{Orbifolds in Mathematics and Physics}, Contemp. Math.,
Amer. Math. Soc., Providence, RI.
{310}, (2002),
135--162.

\bibitem
[Ka3]{orb}
R.~M.~Kaufmann
 {\it Orbifolding Frobenius algebras}.
Int. J. of Math. 14 (2003), 573-619.

\bibitem
[Ka4]   {sq}
R.\ M.\ Kaufmann,  {\it Second quantized Frobenius algebras},
Commun.\ Math.\ Phys.\ to appear.

\bibitem
[Ka5]   {hilb}
R.\ M.\ Kaufmann.
{\it Discrete torsion, symmetric products
and the Hilbert scheme.}
To appear in the conference proceedings of
the conference in honor of Yuri Ivanvovich Manin's 65th birtday.

\bibitem
[Ka6]   {disc}
R.\ M.\ Kaufmann,  {\it The algebra of discrete torsion},
Preprint. math.AG/0208081.

\bibitem[M]{maninbook}
Yu.~Manin, {\it Frobenius manifolds, quantum cohomology, and
moduli spaces}.  Colloquium Publ., {47},  Amer. Math. Soc.,
Providence, RI, 1999.
\bibitem
 [Mo]  {Mo}
S.\ Montgomery.
{\it  Hopf algebras and their actions on rings}.
CBMS Regional Conference Series in Mathematics, 82.
American Mathematical Society, Providence, RI, 1993.

\bibitem
 [P]  {P}
 A.\ Polishchuk.{\it
 Witten's top Chern class on the moduli space of higher spin curves.}
Preprint math.AG/0208112.

\bibitem
 [PV]  {PV} A.\ Polishchuk and A.\ Vaintrob.
{\it Algebraic construction of Witten's top Chern class.}
Advances in algebraic geometry motivated by physics (Lowell, MA, 2000),
229--249, Contemp. Math., 276, Amer. Math. Soc., Providence, RI, 2001.

\bibitem
[S] {Saito} K.\ Saito. {\it Period mapping associated to a primitive
form.} Publ.\ Res.\ Inst.\ Math.\ Sci.\ Kyoto Univ., 19 (1983),
1131--1264. And\\
K.\ Saito.
 {\it Primitive forms for a universal unfolding of a
function with an isolated critical point.} Journ.\ Fac.\ Sci.\
Univ.\ Tokyo, Sec.\ IA. 28 (1982), 775--792.

\bibitem
 [St]  {strachan}
I.A.B.\ Strachan.
{\it Frobenius submanifolds.} J. Geom. Phys. 38 (2001), no. 3-4, 285--307.
And\\
I.A.B.\ Strachan.
{\it Frobenius manifolds and bi-Hamiltonian
structures on discriminant hypersurfaces.}
Integrable systems, topology, and physics (Tokyo, 2000), 251--265,
Contemp. Math., 309, Amer. Math. Soc., Providence, RI, 2002

\bibitem
[V]   {V}
C.\ Vafa,
 String vacua and orbifoldized LG models,
{\it Modern Phys.\ Lett.\ A}{\bf 4} (1989),  1169--1185.


\bibitem
[Wa]   {Wa}
C.\ T.\ C.\ Wall,
 A note on symmetry of singularities,
{\it Bull.\ London Math.\ Soc.}{\bf  12}

\bibitem
[W]{wittenspin}
 E.\ Witten {\it  Algebraic geometry associated with matrix models of two-
dimensional gravity}. Topological models in modern mathematics (Stony Brook,
NY, 1991), Publish or Perish, Houston, TX (1993), 235--269.
and {\it The N-matrix model and gauged WZW models}.  Nucl. Phys. B371
(1992), 191-- 245.1980), 169--175.



\bibitem
[Z]   {zuber}
J.-B.\ Zuber.
{\it
Generalized Dynkin diagrams and root systems and their folding}.
In: Topological field theory, primitive forms and
 related topics (Kyoto, 1996), 453--493, Progr. Math., 160, Birkhäuser Boston, Boston, MA, 1998.

\end{thebibliography}
\end{document}